\documentclass[12pt]{article}
\usepackage{latexsym,amsmath,amssymb,graphics,amscd}

\newcommand{\bfi}{\bfseries\itshape}

\makeatletter \@addtoreset{figure}{section}
\def\thefigure{\thesection.\@arabic\c@figure}
\def\fps@figure{h,t}
\@addtoreset{table}{bsection}

\def\thetable{\thesection.\@arabic\c@table}
\def\fps@table{h, t}
\@addtoreset{equation}{section}

\makeatother

\newtheorem{theorem}{Theorem}[section]
\newtheorem{definition}[theorem]{Definition}
\newtheorem{lemma}[theorem]{Lemma}

\newtheorem{corollary}[theorem]{Corollary}

\begin{document}
\title{Desingularization of Implicit Analytic Differential Equations}
\author{Hern\'{a}n Cendra $^{a}$, Mar\'{\i}a Etchechoury $^{b,1}$ and Alberto Ibort $^{c}$}
\footnotetext[1]{Corresponding author. Fax: $+54-221-424-5875$}
\date{$a$
Universidad Nacional del Sur, Av. Alem 1253\\
8000 Bah\'{\i}a Blanca and CONICET, Argentina.\\
uscendra@criba.edu.ar\\
$b$ Laboratorio de Electr\'onica Industrial, Control e
Instrumentaci\'on,\\
Facultad de Ingenier\'{\i}a, Universidad Nacional de La Plata
and\\
Departamento de Matem\'atica, Facultad de Ciencias Exactas,
Universidad Nacional de La Plata.\\
CC 172, 1900 La Plata, Argentina.\\
marila@mate.unlp.edu.ar\\
$c$ Departamento de Matem\'aticas, Universidad Carlos III de Madrid,\\
Av. de la Universidad 30, Legan\'es, Madrid, Spain.\\
albertoi@math.uc3m.es}

\maketitle


\begin{abstract}
\noindent
The question of finding solutions to given \textsl{implicit
differential equations} (IDE) has been answered by several
authors in the last few years, using different approaches, in an
algebraic and also a geometric setting. Many of those results
assume in one way or another that the subimmersion theorem can be
applied at several stages of the \textsl{reduction algorithm},
which, roughly speaking, allows to reduce a given IDE to a collection
of ODE depending on parameters. The main purpose of the present
paper is to improve some of the known results by introducing at
each stage of the reduction algorithm a \textsl{desingularization}
of the manifolds with singularities that may appear when the subimmersion
theorem cannot be applied. This can be done for analytic IDE by using some fundamental
results on subanalytic subsets and desingularization of closed subanalytic subsets due mainly
to Lojasiewicz, Hironaka, Gabrielov, Hardt, Bierstone, Milman and Sussmann, among others.
We will show how this approach helps to understand the dynamics given by the
Lagrange-D'Alembert-Poincar\'{e} equations for the
\textsl{ symmetric elastic sphere}. 
\end{abstract}

\section{Introduction}\label{INTR}
Implicit differential equations $\phi (x, \dot{x}) = 0$ (IDE) are
very common in science and technology. The case of an ODE $\dot{x}
= f(x)$ is the simplest particular case. The Euler-Lagrange
equations for a given Lagrangian, degenerate or not ~\cite{FM, MS,
CEMARA1}; Lagrange-D'Alembert equations for a nonholonomic system
and their reduced versions Lagrange-D'Alembert-Poincar\'{e}
equations ~\cite{BKMM, CEMARA2, CORTES}; some equations in the
Dirac theory of constraints ~\cite{DIRAC}; electrical circuits
called descriptor or semistate systems ~\cite{new:the}; electric
power systems ~\cite{zab:pow}; nonstandard singularly perturbed
systems ~\cite{McC:sps}; reactive columns ~\cite{kum:dis};
constrained robot systems
~\cite{kri:robot}, are some of the examples.\\

\noindent In the literature, IDE are often called
\textsl{differential algebraic systems} (DAS) in the case in which
the subimmersion theorem can be applied, in one way or another, at
each stage of the so called \textsl{reduction algorithm},
~\cite{rhe2:Impdif}. In the last few years this approach to
finding the solutions to a given DAS in the category of
$C^{\infty}$ manifolds and maps,
and even some more general IDE, has been intensively studied.
Invariants like the so called \textsl{index} of the system were
introduced and calculated, as far as it is possible, for a general
DAS and geometrically motivated algorithms have been developed, see,
for instance, ~\cite{rhe2:Impdif, rhe1:difalg, re1:onag, re2:onae,
rab:gen, sza:geo} just to mention part of the relevant
bibliography more or less connected to the present work.
This kind of approach has now reached a nearly optimal answer to
certain basic questions, for instance, the question about existence
and uniqueness of solutions for DAS. \\

\noindent One may say that the reduction algorithm, as exposed in
the references cited above, uses some of the old ideas of the
Dirac theory of constraints ~\cite{DIRAC}, but generalizes them to
make them useful in a variety of fields beyond Hamiltonian and
quantum mechanics, like control theory and nonholonomic mechanics,
as it appears, for instance,
in the references cited above.\\

\noindent To the best of our knowledge, the singular cases where
the subimmersion theorem cannot be applied have not been
systematically and fully studied in the existing literature,
although there are several related papers.
For instance, an algebraic approach relying
on complex algebraic geometry is presented in
~\cite{pritch:onimp}. In that article the singular cases where the
subimmersion theorem cannot be applied are treated in the case of
implicit differential equations given by complex polynomial
relations. Also an implementation
using computer algebra systems is provided.
Besides ~\cite{pritch:onimp}, cited above, there are
several papers studying specific questions about singular points
and impasse points, see, for instance, ~\cite{rei:sin, soto:imp,
guz:imp, rhe3:imp}. See also ~\cite{BKMM} for examples in
nonholonomic mechanics where a simple example of
desingularization appears.\\

\noindent Desingularization ideas appear in different contexts in
mathematics.
The fundamental theorem of Hironaka ~\cite{hiron1:desing}, whose
proof was simplified and computationally implemented in subsequent
works ~\cite{hau:Hironaka, villa:new}, is a good example of a
desingularization procedure, in this case desingularization of
certain algebraic varieties.
Another fundamental theorem on desingularization was proven by
Bierstone and Milman
~\cite{BIMI}, which actually includes Hironaka's theorem.
See also ~\cite{BM2003}.\\

\noindent In the present work we are going to use some of the
known desingularization results, and also several results from the
theory of subanalytic sets. Our main result shows how to
reduce a given \textsl{real analytic IDE} to a \textsl{real
analytic IDE of constant rank}, defined in section
\ref{sec-difalgsys}, which is considered the simplest case in this
paper.
We do this by conveniently modifying the usual reduction
algorithm by including a desingularization of all the singular
manifolds that may appear at each stage of the algorithm.
The only main result on desingularization that we use for doing this is the theorem of Hironaka on desingularization of
closed subanalytic sets ~\cite{hiron2:subanal}. In fact, the
theorem on desingularization of closed analytic subsets, theorem 5.1 of
Bierstone and Milman in ~\cite{bier88:anal}, is enough for our
purposes. We will keep the conceptual framework
as basic as possible throughout the paper, avoiding the unnecessary usage of schemes.\\

\noindent
Control systems in the category of subanalytic sets have been
recently studied in ~\cite{suss98:control}, where
desingularization techniques, exposed, for instance, in
~\cite{bier88:anal} and ~\cite{suss90:desing},  have been used. In
the present work we also work, in a sense, with control systems in the category of subanalytic sets,
but from the point of view of IDE, which in a sense is dual to the
point of view
of control theory.
In fact, a control system is, roughly, a vector field depending on a parameter, or equivalently, a family of vector fields, or, more generally, a family of local vector fields. On the other hand, as we will see, an IDE of constant rank gives also a family of vector fields, but defined implicitly.
\\

\noindent Working with IDE in the subanalytic category rather than in some
of the $C^{k},$ $k = 1, 2...\infty$ categories, of course is a
limitation, but there are interesting
examples that can be studied in this context. In fact, in the
study of IDE representing several important examples from
mechanics, control theory and other fields, as we have mentioned
before, equations are often given by real analytic
functions.\\

\noindent The relevance of the desingularization method that we
propose to help finding the solutions to a given IDE depends, in
part, on a sufficient knowledge of the geometry of the resulting
desingularizing manifold. For instance, if in a given example the final system obtained by the methods
of the present paper is simply a vector field, existence of equilibrium points may be directly related
to
the topology of the desingularizing manifold.\\

\noindent We develop an interesting example showing how the
desingularization method helps to solve a mechanical system,
namely, the rolling \textsl{symmetric elastic sphere}.
By definition, this system has an \textsl{extra nonholonomic condition} besides
the usual nonsliding condition for a rigid rolling sphere with
only one point of contact with the plane.
This extra condition says that the vertical component of the
angular velocity should be 0.
A physical situation corresponding approximately to this model is
that of an elastic sphere, like, for instance, a rubber sphere,
which is slightly deformed as it rolls on the floor. We assume that
the deformation of the sphere takes place only on that part of the sphere
in contact with the floor.
Then there is an
\textsl{area} of contact, rather than a \textsl{point} of contact, which
imposes, because of the nonsliding condition, the condition that the vertical component of the angular
velocity is 0. From the point of view of reduction theory, this condition introduces singularities in the reduced system of equations, since the standard \textsl{dimension
assumption} for the Lagrange-D'Alembert-Poincar\'{e} equations is
only partially satisfied, see ~\cite{CEMARA2}.
The nonholonomic condition that the vertical component of the angular velocity is 0 has been considered for the example of a rigid body with a fixed point in
~\cite{vesel:newcases}.
We must say that the model that we are considering for a rolling symmetric elastic sphere is not necessarily very realistic because it does not takes into account, for instance, the deformation that may occur in the material far from the area of contact,
or the viscoelasticity, which introduces internal friction. Those are delicate problems in elasticity theory which have been the subject of several investigations for more than a century and are an active area of research in our days, but this is not the subject of our paper. We will concentrate only on the dynamical equations for a rolling sphere under the specified nonholonomic conditions. This kind of system has potential applications
to robotics.
The desingularization process can be performed in detail
in the case of the symmetric elastic sphere and the
desingularizing manifolds can be identified, which leads directly to
solving the system by quadratures.
In fact, we show that the original system is equivalent to a
system on the nonsingular manifold $S^2 \times S^1.$
In spite of the vast literature on the subject of
nonholonomic systems with rolling constraints we believe that our
results in the specific case of the symmetric elastic
sphere are new and interesting by themselves.\\

\noindent
We will work mainly in the category of subanalytic sets
and maps. Manifolds will be usually real analytic manifolds and maps from one manifold to another
will be real analytic maps,
although some of the statement are also valid in a more general context.\\

\noindent In section \ref{sec-difalgsys} we explain some basic
facts about IDE. In section \ref{sec-usingdesing} we describe our
algorithm. In section \ref{sectionmaintheorem} we prove our main
results. In section \ref{section-ERS} we give a detailed
description of the example of the symmetric elastic
sphere.
\section{Implicit differential equations}\label{sec-difalgsys}
\paragraph{Basic notation.}
In this section manifolds will be smooth manifolds and maps from one manifold to another will be smooth maps.
Let $M$ be a given manifold of dimension $n$ and $F$ a vector
space of dimension $m$. Let $a:TM \rightarrow F$ be a map
such that, for each $(x, \dot{x}) \in TM,$ $a(x, \dot{x}) \equiv
a(x)\dot{x} \equiv a(x)(x, \dot{x})$ is linear in $\dot{x}.$ Let $f:M \rightarrow F$ be a
given map.
We will study IDE of the type
\begin{align}
\label{das} a(x)\dot{x}=f(x)
\end{align}\\

\noindent By introducing the trivial vector bundle $M\times F$ we
can think of $a$ as representing a vector bundle map
$$
a : TM \rightarrow M\times F
$$
and of $f$ as being a section of $M\times F.$
Then, for each $x \in M,$ $a(x)$ is a linear map depending
smoothly on $x$ from the tangent space $T_x M$ into the fiber $(x,
F)$
of the trivial bundle.
More generally, we may consider a general vector bundle
with base $M,$ say $\pi : F \rightarrow M$ and an IDE like (\ref{das})
where now $a : TM \rightarrow F$ is a vector bundle map and $f$ is
a section of $F.$ This kind of generalization is important to
describe a sufficiently wide class of IDE.
However, in the present paper we shall describe only the trivial
bundle case, for simplicity and also because it already contains
the essential facts.
The general case can be treated in an essentially similar way.
In this paper the manifold $M$ is called the \textsl{domain} and
the space, or more generally, vector bundle $F$,
is called the \textsl{range} of the IDE.\\

\noindent Given an IDE one has immediately a \textsl{linear algebraic
system} (LAS) for each $x \in M,$ depending smoothly on $x,$ where
the unknown is $(x, \dot{x}),$ for each $x \in M.$ We will call it
the \textsl{LAS associated} to the given IDE.
\\

\paragraph{IDE of constant rank.} Assume that the LAS
associated to (\ref{das}) has solution $(x,\dot{x})$ for each $x
\in M.$ Then (\ref{das}) defines an affine
distribution, generally singular, on $M.$ If, in addition,
$\operatorname{rank} a(x) = \operatorname{rank} [a(x), f(x)]$ is
locally constant, that is, it is constant on each connected
component of $M,$ then the IDE is called an \textsl{IDE of constant
rank}. This is equivalent to saying that the corresponding affine
distribution has \textsl{constant rank} on each connected
component of $M.$
For instance, if $m=n$ and $a(x)$ is invertible for all $x \in M$ then
(\ref{das}) is equivalent to an ODE
\begin{align}
\dot{x} = a(x)^{-1}f(x)
\end{align}
and the rank of the affine distribution is $0$ in this case.\\

\noindent
The case of an IDE of constant rank is the simplest case in our
context and our main result, in section \ref{sectionmaintheorem},
shows that a given \textsl{analytic} IDE can be reduced to a
finite collection of \textsl{analytic} IDE of constant rank, which can be considered also as a single IDE of constant rank in an obvious way.

\paragraph{Reduction of a general IDE to an IDE of the type (\ref{das}).}
It is easy to see that, from the point of view of a general theory where the dimension
does not plays an essential role,
IDE of the type
\begin{align*}
\phi(x,\dot{x})=0
\end{align*}
where the map
$\phi : TM \rightarrow F$ may be nonlinear in $\dot{x}$
are not more general than (\ref{das}).\\

\noindent In fact, let us assume first, for simplicity, that $M$ is an open
subset of a finite dimensional vector space $E.$ An IDE of the
type
\begin{align*}
\phi(x,\dot{x})=0
\end{align*}
can be rewritten in the form
(\ref{das}) with domain $M\times E$
and range
$E\times F$
as follows
\begin{eqnarray*}
 \dot{x} & = & u \\
0 & = & \phi(x,u),
\end{eqnarray*}
which has the form (\ref{das}) with
\[
a(x,u)=\left[\begin{array}{cc} I & 0 \\ 0 & 0 \end{array}\right],
\,\,\, f(x,u)=\left[\begin{array}{c} u \\ \phi(x,u)
\end{array}\right].
\]
The case of a general manifold $M$ and a fiber preserving map
$\phi : TM \rightarrow F$ where $F$ is a vector bundle with base
$M$ can be also reduced to the form (\ref{das}) by an essentially
similar procedure, using a geometric construction involving
pull-backs of bundles.\\

\noindent {\bf Remark.}
(i)
\noindent
We must remark that working with IDE written
in the form (\ref{das}) is an important ingredient of our algorithm, in part
because this form is preserved and the space $F$ remains the \textsl{same}
(or the vector bundle {\sl is replaced} by a pull-back vector bundle)
at each stage of the algorithm,
which simplifies matters as will become evident later.\\
\noindent
(ii)
The system (\ref{das}) can be written
equivalently as follows
\begin{align*}
a(x)\dot{x}=\dot{t}f(x) \\
\dot{t}=1,
\end{align*}
where $t=t(s).$ The first equation of this system can be written
equivalently as follows
\begin{align}
\label{eqhom} b(y)\dot{y}=0
\end{align}
where $y=(x,t)$, $b(y) \in L(T(M\times\Bbb R),\,F),$ $b(y)=[a(x), -f(x)].$\\

\noindent Of course (\ref{eqhom}) is an IDE whose associated LAS
has a solution for each $y,$ but it does not seem that questions
like reachability for (\ref{das}) could be easily reduced to
easily solvable corresponding questions for (\ref{eqhom}). In
other words, this kind of transformation of the system does not
necessarily really simplifies the problems related to a given IDE.
On the other hand, systems like  (\ref{eqhom}) are interesting by themselves
and are related to Pfaffian systems ~\cite{MTA}.

\paragraph{Some notation and operations with IDE.}
It will be convenient to denote $(a, f)$ the IDE
(\ref{das}), from now on. Let $(a, f)$ be a given IDE, say
\[a : TM \rightarrow F,\,\,\,
f : M \rightarrow F.\] Let $N\subseteq M$ be a given submanifold.
The \textsl{restriction} $(a, f)\vert N,$ also written $(a\vert N,
f\vert N),$ of $(a, f)$ to $N$ is defined naturally by the
conditions
$(a\vert N)(x)(x, \dot{x}) = a(x)(x, \dot{x})$
and
$(f\vert N)(x) = f(x),$ for all
$(x, \dot{x}) \in T N.$\\

\noindent Let us remark that the notion of restriction $(a,
f)\vert N$ makes sense also in the case where $N$ is any subset of
$M.$ In fact, we only need to give a meaning to the notion of a
tangent vector $(x_0, \dot{x})$ at a point $x_0 \in N.$ It is simply
the tangent vector in $T_{x_0}M$ to a smooth curve $x(t),$ $t \in
(-\delta, \delta)$ in $M$ such that $x(t) \in N,$
for all $t \in (-\delta, \delta)$ and $x(0) = x_0.$\\

\noindent Let $\varphi : N \rightarrow M$ be a given map.
Then the \textsl{pull-back} $\varphi^{\ast}(a, f) =
(\varphi^{\ast}a, \varphi^{\ast}f)$ is the IDE with domain $N$ and
range $F$ defined by $\varphi^{\ast}a(y)(y,\dot{y}) =
a\left(\varphi(y)\right)\left(T\varphi(y, \dot{y})\right)$ and
$\varphi^{\ast}f (y) = f\left(\varphi(y)\right)$\\

\noindent If $g : F \rightarrow G$ is a linear map we define the
\textsl{projection} of $(a, f)$ \textsl{by} $g$ as being the IDE with domain $M$
and range
$G$
defined by $(g\circ a, g\circ f).$ More precisely, we have
$(g\circ a)(x)(x, \dot{x}) = g(a(x)(x, \dot{x}))$ and $(g\circ
f)(x) = g(f(x)),$ for all
$(x, \dot{x}) \in TM.$\\

\noindent One can define operations like the direct sum
$\oplus$ or tensor product $\otimes$ of IDE in a natural way. For
instance, if $(a_i, f_i),$ $i = 1,2$ are given IDE with domain $M$ and
range $F_i,$ $i = 1,2,$ then we can define the direct sum
$(a_1, f_1) \oplus (a_2, f_2) \equiv (a_1 \oplus a_2, f_1 \oplus f_2)$
as an IDE with domain
$M$
and range
$F_1 \oplus F_2$
by
$(a_1\oplus a_2)(x)
\dot{x} = a_1(x)\dot{x}\oplus a_2(x)\dot{x},$ and
$(f_1 \oplus f_2)(x) = f_1(x)\oplus f_2(x),$ for all
$(x, \dot{x}) \in TM.$
The tensor product is also defined in a natural way.\\

\noindent Using operations like the ones described above, one may
sometimes simplify a given IDE. For instance, working in
coordinates in $F,$ say $(y_1, y_2,...,y_m),$ if some of the
equations, say corresponding to $y_1,$ is a linear consequence of
the others then it can be eliminated by using a projection $g(y_1,
y_2,...,y_m) = (y_2,...,y_m),$ and the resulting system will be
equivalent to the given one. We will need the
following result, whose proof is not difficult.

\begin{theorem}\label{thm0}
Let $(a, f)$ be a given IDE with domain $M$ and range $F$
and let $N\subseteq M$ be a given submanifold defined regularly
by equations
$\varphi = 0,$
where
$\varphi : M \rightarrow H$ and $H$ is a finite
dimensional vector space.
Then the restriction $(a, f) \vert N$ has the same solutions as
the IDE $(a \oplus 0, f\oplus \varphi)$ with domain $M$
and range $F\oplus H.$ It also has the same solutions as the
IDE $(a \oplus D\varphi \oplus 0, f\oplus 0 \oplus \varphi)$
with domain
$M$
and range
$F \oplus H\oplus H.$
Here
$D\varphi : TM \rightarrow H$
is defined by
$D\varphi = p_2 \circ T\varphi,$
where
$p_2 : H \times H \rightarrow H$
is the projection on the second factor
and
$T H \equiv H\times H.$
\end{theorem}

\noindent This theorem is simple but useful. For instance, it
allows sometimes to replace a given IDE by an equivalent IDE whose domain and range are vector spaces, which sometimes simplifies practical calculations avoiding
the usage of local charts whenever it is convenient.
More precisely, let a given IDE $(a, f)$ having domain
$M\subseteq L,$ imbedded in the vector space $L$ and defined
regularly by an equation $\varphi = 0,$ where $\varphi : L
\rightarrow H,$ and let $a$ be defined by a restriction $a =
A\vert TM,$ where $A : TL \rightarrow F.$ Then, according to
theorem \ref{thm0}, one can work equivalently with the system
$(A\oplus 0, f\oplus \varphi),$
whose domain and range are vector spaces.\\

\noindent One can obviously define a category whose objects are
of the type
$\left(M, F, (a, f)\right),$ where $M$ is a manifold, $F$ is a
vector bundle over $M$ and $(a, f)$ is an IDE with domain $M$ and
range $F.$ A morphism $\varphi : \left(M, F, (a, f)\right)
\rightarrow \left(N, G, (b, g)\right)$ is given by a map
$\varphi_d : M \rightarrow N$ and a vector bundle map $\varphi_r :
F \rightarrow G$
over
$\varphi_d$
such that, for any $(x, \dot{x}) \in M,$ we have
$\varphi_r \left(a(x)\dot{x}\right) = b(y)\dot{y},$ and $g(y) =
\varphi_r\left(f(x)\right),$ where $(y, \dot{y}) = T\varphi_d (x,
\dot{x}).$ However,
although this perspective is interesting,
we will not use the categorical context in
this paper since, as we have said before, we want to keep the context as basic as possible.

\paragraph { Basic reduction algorithm for solving an IDE.}

We have the following basic algorithm to solve  (\ref{das})
recursively. It is a reformulation in our context of essentially
known ideas contained in the references given before.
However, we must remark that the fact that we have chosen to write
a given IDE in the form (\ref{das}) where the space (or more
generally, vector bundle) $F,$ remains unchanged
(or is replaced by a pull-back bundle)
throughout the
reduction process, briefly described next, presents some clear technical advantages.
\\

$(a_1)$ $M_1 \subseteq M$ is the subset of all $x \in M$ such that
the
LAS (\ref{das}) has solution.\\

$(a_2)$ For $k=1,2,\ldots$ we assume that $M_k$ is a submanifold
and $M_{k+1} \subseteq M_{k}$ is the subset of all $x \in M_k$
such that the LAS (\ref{das}) has a solution
$(x,\dot{x}) \in TM_{k}.$\\

\noindent If in a given example the assumptions that $M_{k}$ is a
submanifold made at each stage of the previous algorithm are
satisfied then the algorithm itself stabilizes at a certain stage
$q$, that is $M_{q}= M_{q+1}.$ Then the system (\ref{das})
restricted to $M_{q}$ has a solution $(x,\dot{x}) \in TM_{q}$, as
a LAS, for each $x \in M_{q}.$ Thus we obtain an IDE with domain $M_{q}$
and range
$F$
which may or may not be of constant rank according to the
definition given before, but it is equivalent to the original one and
the associated LAS has solution for all
$x \in M_q.$
The assumption that $M_{k + 1} \subseteq M_k$ is a
submanifold is obtained in practice by using, somehow, the
subimmersion theorem or equivalent results.
An important case in which this scheme works is the case where $M$
is a vector space, $a(x)$ is independent of $x$ and $f(x)$ is
linear in $x.$ This kind of situation was studied in
~\cite{rabier:lineal}. It can be also studied with our methods,
leading to perfectly identifiable and simple desingularizing
systems, which we are not going to explain in the present work.
\\
\section{Desingularization}\label{sec-usingdesing}
The hypothesis that $M_k$ is a submanifold at each stage of the
basic algorithm described above is too restrictive as it is not satisfied for many examples of interest.
Some rather general geometric approaches to IDE in the $C^{\infty}$
category, like in ~\cite{rhe2:Impdif, sza:geo}, also assume some
kind of geometric version of the subimmersion theorem.
We are going to show that in order to overcome part of those
limitations, and at the cost of working in the subanalytic
category rather than the $C^{\infty}$ category, one can use
results from real analytic desingularization theory.
This is possible thanks to the now well established theory of
semianalytic and subanaliytic sets developed originally by
Lojasiewicz ~\cite{lojasiewicz1, lojasiewicz2, lojasiewicz3}.
Important results in this field and systematic expositions using
techniques which are simpler than the original ones are due to
Gabrielov ~\cite{Gabrielov}, Hironaka  ~\cite{hiron2:subanal},
Hardt ~\cite{Hardt1, Hardt2}, Bierstone and Milman
~\cite{bier88:anal}, Sussmann ~\cite{suss90:desing}, and others.
Our main reference will be ~\cite{bier88:anal}, where an excellent
and very readable exposition of important points of the theory of
subanalytic sets has been written.
\\

\noindent
For the rest of this paper manifolds and maps will be real analytic, unless otherwise specified. For instance,
if
$(a, f)$
is a given IDE with domain
$M$
then
$M$
will be a real analytic manifold and
$a,$
$f$
will be real analytic maps.
\begin{definition}
Let $M$ be a real analytic manifold and let $X$ be a closed
subanalytic subset of $M.$ A desingularization of $X$ is a a real
proper analytic map $f : N \rightarrow M$ such that $f(N) = X,$
where $N$ is a real analytic manifold of the same dimension as
$X.$
\end{definition}
This is a relatively weak notion of desingularization, but it is
enough for our purposes. Existence of desingularizations $f: N
\rightarrow M$ is guaranteed by the following theorem of Hironaka,
see ~\cite{bier88:anal, suss90:desing} and references therein.
\begin{theorem}\label{theoremH1}
Let $M$ be a real analytic manifold and let $X$ be a closed
subanalytic subset. Then there is a desingularization $f : N
\rightarrow M$ of $X.$
\end{theorem}
Desingularization results that include those of Hironaka have been recently proved in Bierstone and Milman ~\cite{BIMI}.\\

\noindent In fact, in the present paper we only need the following weaker
desingularization result, which is theorem 5.1 in
~\cite{bier88:anal},
\begin{theorem}\label{theoremH23}
Let $M$ be a real analytic manifold and let $X$ be a closed
analytic subset. Then there is a desingularization $f : N
\rightarrow M$ of $X.$
\end{theorem}
\paragraph{Description of the algorithm.}
Let $M$ be a manifold of dimension $d$ and let $(a,
f)$ be a given IDE with domain $M$ and range $F.$
The basic result proved in this paper to solve the IDE (\ref{das})
consists, roughly, in transforming it into an equivalent IDE, say
\begin{align*}
\tilde{a}_2(y)\dot{y}=\tilde{f}_2(y)
\end{align*}
on a manifold $\tilde{M}_2$, which is a DAS of constant rank.
The manifold $\tilde{M}_2$ will be constructed by an algorithm that
involves
a desingularization process.\\

\noindent
{\bfi The decomposition $M = M_0 \cup M_1 \cup M_2.$}
First, let us assume that $M$ is a connected manifold of dimension
$d$. For $i = 0, 1, \ldots,$ let
\begin{align*}
S_i(M) &=
\{x \in M\,|\,\operatorname{rank}a(x) \leq i\}\\
&= \{x \in M\,|\,\operatorname{det}A(x)=0, \,\, \mbox{$A(x)$
submatrix of $a(x)$ of order $i+1$}\}
\end{align*}
$S_i(M)$ is clearly a closed analytic subset of $M$, defined by
analytic
equations, for $i=0, 1,\ldots.$\\

\noindent
Also, for $i = 0, 1,...,$ let $L_i(M) \subseteq S_i(M)$ be defined
by
\begin{align*}
L_i(M) &=
\{x \in S_i(M)\,|\,\operatorname{rank}[a(x),f(x)] \leq i\}\\
&= \{x \in S_i(M)\,|\,\operatorname{det}A(x)=0, \mbox{$A(x)$
submatrix of $[a(x),f(x)]$
 of order $i+1$}\}.
\end{align*}
Each $L_i(M)$ is a closed analytic subset of $M$ defined by analytic equations.\\

\noindent
Let
\[
S_{k_1}(M) \subset S_{k_2}(M) \subset \ldots \subset S_{k_{r}}(M),
\]
where
$k_r = k_r(M),$
be the distinct nonempty $S_i(M).$
We observe that $S_{k_r}(M) \equiv M.$ Consider the corresponding
inclusions
\[
L_{k_1}(M) \subseteq L_{k_2}(M) \subseteq \ldots \subseteq
L_{k_{r}}(M).
\]
We have that
$\operatorname{rank}a(x)=\operatorname{rank}[a(x),f(x)]=k_j$ for
each $x \in L_{k_j}(M)-S_{k_{j-1}}(M)$, $j=1,\ldots,r.$ The LAS
associated to (\ref{das}) has solution for each $x \in
L_{k_j}(M)-S_{k_{j-1}}(M)$, $j=1,\ldots,r,$ where we have, by definition,
$S_{k_0} = \emptyset.$\\

\noindent We remark the following useful facts: the set
$L_{k_j}(M)-S_{k_{j}-1}(M)$ may be empty, for some $j=1,\ldots,r;$
we have $\operatorname{dim}S_{k_{r-1}}(M) < \operatorname{dim} M;$
if $\operatorname{dim}(L_{k_r}(M)) = d,$ then $L_{k_r}(M) = M.$
\\

\noindent Now let $M$ be a manifold of dimension $d$ and assume
that
$$
M_m = \bigcup_{j} W_j
$$
is the union of the connected components of $M$ of maximal
dimension
$d.$\\

\noindent We will consider the following pairwise disjoint
conditions for a given $W_j \subseteq M_m,$
\begin{description}
\item (a) $L_{k_{r}}(W_j) = \emptyset$ \item (b) $L_{k_{r}}(W_j)
\not= \emptyset$ and $\operatorname{dim}L_{k_{r}}(W_j) < d$ \item
(c) $L_{k_{r}}(W_j)  \not= \emptyset$ and
$\operatorname{dim}L_{k_{r}}(W_j) = d.$
\end{description}
We now define the following pairwise disjoint subsets of $M.$
\begin{align*}
M_0 &= (M - M_m) \cup \bigcup_{b} L_{k_r}(W_j) \cup
\bigcup_{c}S_{k_{r-1}}(W_j)\\
M_1 &= \bigcup_{a} W_j \cup
\bigcup_{b} \left(W_j - L_{k_r}(W_j)\right)\\
M_2 &= \bigcup_{c} \left(W_j - S_{k_{r-1}}(W_j)\right).
\end{align*}
We have the following assertions, whose proof is easy:
each subset $L_{k_r}(W_j) \subseteq W_j,$ and
each subset
$S_{k_{r-1}}(W_j)
\subseteq W_j,$ is a closed analytic subset of $W_j$ defined by
analytic equations on $W_j.$ In consequence, $W_j - L_{k_r}(W_j),$
$W_j - S_{k_{r-1}}(W_j)$ are open submanifolds of
$W_j.$\\

\noindent The manifold $M$ is the disjoint union
$$
M = M_0 \cup M_1 \cup M_2.
$$
The manifolds $M_1$ and $M_2$ are open submanifolds of $M.$ The
subset $M_0$ is a union of subsets defined by analytic equations
on each $W_j$, union $M - M_m,$ and we have that
$\operatorname{dim}M_0 < d.$\\

\noindent
{\bfi Restrictions $(a, f)\vert M_0,$ $(a, f)\vert M_1$, $(a, f)\vert M_2$ and desingularization of $(a, f)\vert M_0.$}
We have that the LAS associated to (\ref{das}) has no solution for
$x \in M_1.$ On the other hand, it has solution for all $x \in
M_2,$ moreover, $(a, f)\vert M_2,$
\textsl{is an IDE of constant rank}.\\

\noindent It remains to see what happens with the system
restricted to $M_0.$ The idea here is to desingularize each
closed analytic subset $L_{k_r}(W_j) \subseteq W_j,$ and
$S_{k_{r-1}}(W_j) \subseteq W_j.$
By forming the disjoint union of those desingularizations and
$M-M_m$ one obtains a desingularization of $M_0$ say
\[
{\pi}_0:M^1 \rightarrow M,\,\mbox{where}\,\,\pi_0(M^1)=M_0.
\]

\noindent Then (\ref{das}) restricted to $M_0$, that is
$(a,f)|M_0,$ can be naturally \textsl{lifted},
using the pullback operation, to an IDE
$(a_1,f_1)=\pi_0^*\left((a,f)|M_0\right)$ on $M^1$ as follows
\begin{align*}
a_1(y)\dot{y}=a(\pi_0(y))T_y\pi_{0}(y, \dot{y})\\
f_1(y)=f(\pi_{0}(y)).
\end{align*}

\noindent We should remark at this point (see also the paragraph
{\sl Some notation and operations with IDE}, in section \ref{sec-difalgsys})
that in the present paper a tangent
vector $(x, \dot{x})$ to $M_0$ at a point $x \in M_0$ is a vector
$(x, \dot{x}) \in T_x M$ such there is an analytic curve $z(t) \in
M_0,$ say defined for $t \in (-\delta, \delta),$ such that
$z(0) = x$ and the
tangent vector to $z(t)$ at $t = 0$ as a curve in $M$ coincides
with $(x, \dot{x}).$
In particular, if $y(t)$ is a given analytic curve in $M^1$ then
$$
T_y \pi_0 (y, \dot{y}) = \left.\frac{d \pi_0
\left(y(t)\right)}{dt}\right|_{t = 0},
$$
is a tangent vector to $M_0$ at $\pi_0 \left(y(0)\right).$
\\

\noindent Note that $M^1$ is a manifold of dimension
$\operatorname{dim}{M^1}=
\operatorname{dim}{M_0} < d.$\\

\noindent
{\bfi Complete desingularization of $(a, f)$ in a finite
number of steps.} Now we repeat the process for the IDE
$(a_1,f_1)$ with domain ${M^1}$ and range $F,$ proceeding as we did
before with the system $(a, f)$ with domain $M$
and range $F.$
We obtain a
decomposition
\[M^1=M^1_0 \cup M^1_1 \cup M^1_2.\]
We know that there is no solution to the LAS system
\[a_1(y)\dot{y}=f_1(y)\]
for $y \in M_1^1.$
We also know that there is solution to the same LAS
system for $y \in M_2^1,$ moreover, $(a_1, f_1)\vert M^{1}_{2}$ is
an IDE of constant rank. Now we desingularize $M_0^1$
\[\pi_1:M^2 \rightarrow M^1,\,\pi_1(M^2)=M_0^1\]
and repeat the process. Finally, we obtain a finite sequence of
manifolds and maps
\[
M^q \stackrel {\pi_{q-1}} {\rightarrow} M^{q-1} \stackrel
{\pi_{q-2}} {\rightarrow} \ldots \stackrel {\pi_1}{\rightarrow}
M^1 \stackrel {\pi_0}{\rightarrow} M,
\]
where $\pi_0 (M^1) = M_0,$ $\pi_1 (M^2) = M_0^1,$ and in general
$\pi_i (M^{i + 1}) = M_0^i,$ for $i = 0,...,q-1,$ where we have
written $M^0 \equiv M$
to unify the notation.\\

\noindent We have obtained a finite recursive procedure that
reduces the problem to a finite number of IDE of constant rank,
namely, the IDE of constant rank $(a_i,f_i)\vert M_2^i,$ for $i =
0,1,...,q,$ where we have written $(a_0, f_0) = (a, f),$ to unify
the notation. We will call this a \textsl{desingularization
process} and the sequence of maps $\pi_i$ and IDE $(a_{i+1},
f_{i+1}),$ $i
= 0,...,q-1$ a \textsl{desingularization} of $(a, f).$\\

\noindent
The collection $(\tilde{a}_2, \tilde{f}_2)$
of IDE
$(a_k, f_k)\vert M^k_2,$ $k = 0,...,q,$
defines a single IDE of
constant rank in the disjoint union
$\tilde{M}_2 = \bigsqcup_{k = 0}^{q}
M^k_2,$
as we have said at the beginning of the paragraph
\textsl{Description of the algorithm}.
We have a natural projection
$\tilde{\pi}_2 : \tilde{M}_2 \rightarrow M.$
This IDE
$(\tilde{a}_2, \tilde{f}_2)$
with domain
$\tilde{M}_2$
and range
$F$
is called the \textsl{desingularizing IDE}.
\\

\noindent \textbf{Remark.}
As we have said before the range
$F$
remains the same throughout the application of the algorithm. However, in practice it is sometimes convenient to apply theorem \ref{thm0}, which may imply a change of
$F,$
to simplify calculations.
\section{The Main Results}\label{sectionmaintheorem}
In this section we will show in which precise sense the solutions
to the desingularizing system
$(\tilde{a}_2, \tilde{f}_2)$
of a given analytic IDE $(a, f)$ are
related to the solutions to $(a, f).$
It is clear that in certain
simple examples of IDE one can show in a more or less direct way that
solutions to the desingularizing system
project via $\tilde{\pi}_2$ onto
solutions to the given IDE, and also that
solutions to the IDE are such projections.
For instance, if one is interested in the local behavior of solutions near
a singular point
of
$M_0$
it is sometimes enough to use a simple blow-up to desingularize
$M_0$
at that point, and
one can show in certain cases
in a simple and useful way the relationship between solutions to
the given system and solutions to the system obtained by blow-up, see
~\cite{BKMM}.
However, in
this paper we want a more precise and general global result showing that a curve, belonging to a certain convenient class of curves, is
a solution to a given IDE if and only if it is essentially the projection via the map
$\tilde{\pi}_2$
of a solution to the desingularizing system belonging to the same class of curves.
This is important, for instance, if
one is interested in global aspects of solutions, like extension of solutions, or in a description of the family of all solutions.
In order to be able to use
the theory of subanalytic sets we need to define carefully a convenient class of curves, which we will do in the next paragraph.\\

\noindent From now on, we will often use the theorem 6.1
of ~\cite{bier88:anal} that a subanalytic subset of dimension 1 of
an analytic manifold $M$ is a semianalytic subset. We will also
often use the fact that the image of a relatively compact
subanalytic subset under a subanalytic map is a subanalytic subset, see
~\cite{bier88:anal} immediately after definition 3.2.
Using this we can deduce that the image of a relatively compact
subanalytic subset under an analytic map is a subanalytic subset, which will be also useful for us.\\

\noindent In what follows, we will work with several types of intervals, like
$(\tau_0, \tau_1),$ $[\tau_0, \tau_1),$  $(\tau_0, \tau_1]$
or
$[\tau_0, \tau_1].$ We will usually assume that
$\tau_0$ and
$\tau_1$
are real numbers. However, some definitions and results are valid also for the case in which the open end of an interval is
$\pm \infty,$
that is,
for intervals
$(-\infty, \tau_1);$ $(\tau_0, +\infty);$
$(-\infty, +\infty);$ $[\tau_0, +\infty);$ $(-\infty, \tau_1].$
\paragraph{The notion of an as-curve.}
Inspired by ~\cite{bier88:anal}, definition 3.2, we will define
\begin{definition}\label{defsubcurv}
A subanalytic curve
$x : (t_0, t_1) \rightarrow M$
is a subanalytic map, that is, a map such that
$\operatorname{graph}x \subseteq \Bbb{R}\times M$
is a subanalytic subset. We define the notion of a subanalytic curve
$x : [t_0, t_1) \rightarrow M,$ $x : (t_0, t_1] \rightarrow M$
or
$x : [t_0, t_1] \rightarrow M$
in a similar way.
\end{definition}
In definition \ref{defsubcurv} since
$\operatorname{dim} (\operatorname{graph} x) = 1$ we have that
$\operatorname{graph} x$
is a semianalytic set.
\begin{lemma}\label{lemmasuban1}
$(a)$
Let
$x : [t_0, t_1) \rightarrow M$
be a continuous subanalytic curve whose graph is a relatively compact subset.
Then there is a uniquely defined continuous subanalytic extension
$\bar{x} : [t_0, t_1] \rightarrow M.$
A similar result holds for subanalytic curves
$x : (t_0, t_1] \rightarrow M$  or $x : (t_0, t_1) \rightarrow M.$\\
$(b)$
Let
$x : [t_0, t_1) \rightarrow M$
be a continuous subanalytic curve whose graph is not a relatively compact subset.
Then
$\operatorname{graph} x$
is closed.
A similar result holds for subanalytic curves
$x : (t_0, t_1] \rightarrow M$  or $x : (t_0, t_1) \rightarrow M.$
\end{lemma}
\textbf{Proof.}
First we shall prove
$(a).$
We need to show first that the limit of
$x(t)$ as
$t \rightharpoonup t_1^{-}$
exists.
Using corollary 2.8 of ~\cite{bier88:anal} we can deduce that the closure
$\overline{G}$
of
$G = \operatorname{graph}x$
in
$\Bbb{R} \times M$
is subanalytic and compact, and then also
$\overline{G} \cap (\{t_1\} \times M)$
is subanalytic and compact.
Since
$\overline{G} \cap (\{t_1\} \times M)$
is nonempty let
$x_1 \in \overline{G} \cap (\{t_1\} \times M)$
be given.
One can choose local coordinates at
$x_1,$
and for any small
$\epsilon > 0$
the set
$G_{\epsilon} = G \cap \left((t_0, t_1) \times B_{\epsilon}(x_1)\right)$
is relatively compact.
We can easily deduce from theorem
3.14 of ~\cite{bier88:anal}
that any relatively compact subanalytic set has a finite number of connected components.
Let
$C_i,$
$i = 1,...,n(\epsilon)$
be the connected components of
$G_{\epsilon}.$
It is not difficult to see that each connected component
$C_i$
is of the type
$C_i = \operatorname{graph}\left(x \vert (\alpha_i, \beta_i)\right),$
$i = 1,...,n(\epsilon).$
We can assume without loss of generality that
$\beta_i \leq \alpha_{i + 1},$
$i = 1,...,n(\epsilon) - 1.$
Since
$x_1$
is a limit point of
$G_{\epsilon}$
we must have
$\beta_{n(\epsilon)} = t_1,$
which implies that
$x(t) \in B_{\epsilon}(x_1)$
for all
$t \in (\alpha_{n(\epsilon)}, \beta_{n(\epsilon)}),$
as we wanted to prove.
The fact that the continuous extension
$\bar{x} : [t_0, t_1] \rightarrow M$
is a subanalytic curve follows from corollary 2.8 of
~\cite{bier88:anal}.
The rest of the proof of
$(a)$ can be performed in a similar way.
Now we will prove
$(b).$
If
$\operatorname{graph}x$
is not relatively compact then
$\overline{G} \cap (\{t_1\} \times M)$
must be empty, otherwise we can proceed as in the proof of
$(a)$
and we can conclude that
the limit of
$x(t)$
as
$t \rightharpoonup x_1^-$
exists and then, one can show
that
$G$
is relatively compact.
Since all the limit points of
$G$
not in
$G$
must belong to
$\overline{G} \cap (\{t_1\} \times M)$
we have that
$G$
is closed.
The rest of the proof of
$(b)$ can be performed in a similar way.
\quad $\blacksquare$\\

\noindent
In order to define a convenient class of curves to solve a given IDE we introduce the following notion
\begin{definition}
(a) An \textsl{analytic-semianalytic-curve} $x(t),$ $t\in (t_0,
t_1)$, in $M,$ where $M$ is a given manifold, is an
analytic map $x : (t_0, t_1) \rightarrow M$
which is also a subanalytic curve, that is,
such that
$\operatorname{graph}x$ is a semianalytic subset of $\Bbb{R}\times
M.$ We will often call such an analytic-semianalytic-curve in $M$
simply an \textsl{as-curve in $M.$}\\
\noindent $(b)$ An \textsl{analytic-semianalytic-curve} (or
\textsl{as-curve}) $x(t),$ $t\in [t_0, t_1)$, $(t\in (t_0, t_1],
t\in [t_0, t_1])$ in $M,$ where $M$ is a given manifold,
is a continuous map $x : [t_0, t_1) \rightarrow M,$ (respectively,
$x : (t_0, t_1] \rightarrow M,$ $x : [t_0, t_1] \rightarrow M$), which is also a subanalytic curve, that is,
such that $\operatorname{graph}x$ is a semianalytic subset of
$\Bbb{R}\times M,$ and $x\vert (t_0, t_1)$ is an as-curve in $M.$
\end{definition}

\noindent For instance, $x = \sqrt{t}, \, t \in (0, c),$ (or $t
\in [0, c),$ $t \in (0, c],$ $t \in [0, c]$), with $c > 0,$ are
as-curves in $\Bbb{R}$.
On the other hand, $x = t \sin (\pi/
t), \,t \in (0, c),$ with $c > 0,$ is \textsl{not} an as-curve in $\Bbb{R},$ but $x = t \sin (\pi/ t), \,t \in (\delta, c),$
with $0 < \delta < c ,$ \textsl{is} an as-curve in
$\Bbb{R}.$ \\

\noindent Next, we will give some lemmas where some basic properties of as-curves, that we need to prove our main
results, are proved.
\begin{lemma}\label{lemaa0}
$(a)$
Let $x : (t_0, t_1) \rightarrow N$ be a
given analytic map, where $N$ is a given manifold. Then
any map $x\vert (\bar{t}_0, \bar{t}_1): (\bar{t}_0, \bar{t}_1) \rightarrow N,$ where
$\bar{t}_0$ and $\bar{t}_1$ are such that $t_0 < \bar{t}_0
<\bar{t}_1 < t_1,$ is an as-curve in $N.$
In a
similar way, any map $x\vert [\bar{t}_0, \bar{t}_1): [\bar{t}_0, \bar{t}_1)\rightarrow N,$ $x\vert (\bar{t}_0, \bar{t}_1]:
(\bar{t}_0, \bar{t}_1]\rightarrow N,$ or
$x\vert [\bar{t}_0, \bar{t}_1]: [\bar{t}_0,
\bar{t}_1]\rightarrow N,$ with
$\bar{t}_0$ and $\bar{t}_1$ as before is an as-curve in $N.$\\
$(b)$
Let $x : (t_0, t_1) \rightarrow N$
be an as-curve
and assume that
$\operatorname{graph} x$
is a relatively compact subset of
$\Bbb{R} \times N.$
Then there is a unique continuous extension
$\bar{x} : [t_0, t_1] \rightarrow M$
which is an as-curve.
\end{lemma}
\textbf{Proof.} Let us prove
$(a).$
We have that $\operatorname{graph}(x\vert
\left(\bar{t}_0, \bar{t}_1)\right)$ is a semianalytic subset of
$\Bbb{R}\times N$
defined as
$\{(t, x)\in (t_0, t_1)\times N:
\bar{t}_0 < t <\bar{t}_1 \, , \, x = x(t)\}.$ The rest $(a)$ can be proved in a similar way.
To prove $(b)$ we simply apply lemma \ref{lemmasuban1}
\quad $\blacksquare$
\begin{lemma}\label{lemaasa}
Let $x : [t_0, t_1) \rightarrow M$ be a subanalytic map and assume that $x$ is continuous at
$t_0.$ Then there exists $t_2 \in (t_0, t_1)$ such that $x \vert
[t_0, t_2]$ is an as-curve.
A similar result holds for a
subanalytic map $x : (t_0, t_1] \rightarrow M$ continuous at
$t_1,$
that is, there exists $t_2 \in (t_0, t_1)$ such that $x \vert
[t_2, t_1]$ is an as-curve.
\end{lemma}
\textbf{Proof.}
If
$x$
is a constant the result follows immediately. Let us assume that
$x$
is not a constant.
It is easy to see that
$\operatorname{graph} (x\vert [t_0, t_2])$
is a semianalytic subset of
$\Bbb{R} \times M$ of dimension 1,
for every
$t_2 \in (t_0, t_1).$
We can assume without loss of generality, using, for instance, Whitney embedding theorem, that
$M \subseteq U$
is an analytic submanifold of
$U,$
where
$U$
is a real finite dimensional vector space.
Let
$p_1 : \Bbb{R} \times U \rightarrow \Bbb{R}$
be the projection onto the first factor.
Continuity of
$x$
at
$t_0$
implies that one can
choose
$b \in (t_0, t_1),$
such that
$\operatorname{graph} (x\vert [t_0, b])$
is a relatively compact subanalytic subset of
$\Bbb{R} \times U.$
According to lemma 3.4 of ~\cite{bier88:anal} we have that
$\operatorname{graph} (x\vert [t_0, b])$
is a finite union of connected smooth semianalytic subsets
$A$
such that, for each
$A,$
$\operatorname{rank}(p_1 \vert A)$
is constant.
It is not difficult to see that each
$A$
has dimension 0 or 1 and that
$\operatorname{rank}(p_1 \vert A)$
is 0 or 1. Moreover, we can see that there must be an $A,$
say
$A = A_0,$
such that
$\operatorname{rank}(p_1 \vert A_0) = 1,$
$p_1 (A_0) = (t_0, t_2),$ for some
$t_2 \in (t_0, b],$
and
$p_1 (\bar{A}_0) = [t_0, t_2].$ From this we can easily deduce
that
$x\vert (t_0, t_2)$
is an as-curve and moreover, using lema \ref{lemaa0}, $(b),$
that
$x\vert [t_0, t_2]$
is an as-curve. The rest of the proof can be performed in a similar way.
\quad $\blacksquare$

\begin{lemma}\label{lemmaas0}
$(a)$ Let $x : [t_0, t_1]\rightarrow M,$ $x : [t_0, t_1)\rightarrow M$ or $x
: (t_0, t_1]\rightarrow M$ be an as-curve in $M.$
Then $x\vert (t_0, t_1)$ is an as-curve in $M.$\\
\noindent
$(b)$ Let $x : (t_0, t_1)\rightarrow M$ be an as-curve in $M$ and assume that there is a continuous
extension $\bar{x} : [t_0, t_1]\rightarrow M,$ $\bar{x} : [t_0, t_1)\rightarrow M$ or $\bar{x} : (t_0, t_1]\rightarrow M.$ Then
$\bar{x}$ is an as-curve in $M.$\\
\noindent
$(c)$ Let
$f : M \rightarrow N$ be a given analytic map.
Let $x(t),\, t \in
[t_0, t_1)$
$\left(t \in (t_0, t_1], t \in [t_0, t_1],\right.$
$\left.t \in (t_0, t_1)\right)$
be an as-curve in $M.$
Then $f\left(x(t)\right)$, $t \in [\bar{t}_0, \bar{t}_1],$
is an as-curve in $N,$ for each $[\bar{t}_0, \bar{t}_1] \subseteq [t_0, t_1)$,
(respectively, $[\bar{t}_0, \bar{t}_1] \subseteq (t_0, t_1]$, $[\bar{t}_0, \bar{t}_1] \subseteq [t_0, t_1],$ $[\bar{t}_0, \bar{t}_1] \subseteq (t_0, t_1)$).
If
$\operatorname{graph} x$
is relatively compact then
$f \left(x(t)\right),\, t \in
[t_0, t_1),\,\left(t \in (t_0, t_1], t \in [t_0, t_1], t \in (t_0, t_1)\right)$
is an as-curve
and
$\operatorname{graph} (x\circ f)$
is relatively compact.
\end{lemma}
\textbf{Proof.}
Part $(a)$ is an immediate consequence of the definitions. Part $(b)$ follows
easily using corollary 2.8 of ~\cite{bier88:anal}.
To prove $(c)$ observe first that $f \circ x : (t_0, t_1) \rightarrow
N$ is an analytic map.
Since $\operatorname{graph} (x\vert [\bar{t}_0, \bar{t}_1])$ is a
semianalytic compact subset of $\Bbb{R} \times M$
of dimension $1$
we
have that $\operatorname{graph}(f\circ x\vert [\bar{t}_0, \bar{t}_1]) =
\left(1_{\Bbb{R}}\times f\right) \left(\operatorname{graph}(x\vert [\bar{t}_0, \bar{t}_1]\right),$
taking into account that
$1_{\Bbb{R}} \times f : \Bbb{R} \times M \rightarrow \Bbb{R}
\times N$ is an analytic map,
is also semianalytic, and a similar proof can be given for the case of the
intervals
$(t_0, t_1],$
$[t_0, t_1],$
$(t_0, t_1).$
If
$\operatorname{graph}x$
is relatively compact then according to lemma \ref{lemmasuban1}
we have a continuous extension
$\bar{x}$
which is an as-curve
and therefore we can apply the first part of
$(c)$
to this extension
The rest of the proof follows easily.
\quad $\blacksquare$\\

\begin{lemma}\label{lemmaas1}
Let $N$ be a given manifold and let $x :
[t_0, t_1) \rightarrow N$ be an as-curve in $N$ which is not a
constant. Then there exist $t_2 \in (t_0, t_1)$  such that $x\left((t_0, t_2)\right)$ is an analytic
submanifold which is also a semianalytic subset, $x \vert (t_0,
t_2) : (t_0, t_2) \rightarrow x\left((t_0, t_2)\right)$ is an
analytic diffeomorphism and $x \vert [t_0, t_2] : [t_0, t_2]
\rightarrow x\left([t_0, t_1]\right)$ is an homeomorphism. Similar
results hold for an as-curve $x : (t_0, t_1] \rightarrow N.$
\end{lemma}
\textbf{Proof.}
Let
$t_0 < \bar{t}_1 < t_1,$ then we have that
$\operatorname{graph}(x\vert [t_0, \bar{t}_1])$
is a compact semianalytic subset of
$\Bbb{R} \times N$
of
dimension $1.$
Let $x_i,$ $i = 1,...,n$ be local analytic coordinates centered at
$x(t_0)$ therefore $x_i(t_0) = 0$ for $i = 1,..,n.$ Without loss
of generality we can assume that $x_i(t)$ is defined for all $i =
1,..,n$ and all $t \in [t_0, \bar{t}_1].$
For some index, say $j \in
\{1,...,n\},$ we must have that $x_j(t)$ is not a constant.
We are going to show that there exists $t_2 \in (t_0, \bar{t}_1]$ such that
$x_j : [t_0, t_2] \rightarrow \Bbb{R}$ satisfies certain conditions from which the lemma follows.
Since
$x_j(t)$
is the projection on the
$j$-coordinate axis of the curve
$x\vert [t_0, \bar{t}_1]$
using lemma \ref{lemmaas0}, $(c)$
we have that
$\operatorname{graph}\left(x_j\vert [t_0, \bar{t}_1]\right)
\subseteq \Bbb{R} \times \Bbb{R}$ is a
semianalytic compact subset of dimension $1.$
The restriction
to
$\operatorname{graph}(x_j \vert [t_0, \bar{t}_1])$
of the projection
$p_2 : \Bbb{R} \times
\Bbb{R} \rightarrow \Bbb{R}$ onto the second factor satisfies
$p_2 \left(t, x_j(t)\right) =
x_j(t),$ for all $t \in [t_0, \bar{t}_1].$
Then according to lemma 3.4
of ~\cite{bier88:anal}  $\operatorname{graph}\left(x_j \vert [t_0, \bar{t}_1]\right)$
is a finite
union of connected smooth semianalytic subsets $A$ such that
$\operatorname{rank}(p_2 \vert A)$ is constant on $A,$ and it is
easy to show that $\operatorname{rank}(p_2 \vert
A)$ can take only the values $0$ or $1.$
It is also easy to see
that for at least one such $A$ one must have that
$\operatorname{rank}(p_2 \vert A) = 1,$ $p_2(A) = x_j\left((t_0, t_2)\right)$
for some $t_2 \in (t_0, \bar{t}_1]$ and $\overline{p_2(A)}
=
p_2(\overline{A})  = x_j \left([t_0, t_2]\right),$
and therefore that
$x_j((t_0, t_2))$
is an open interval.
Moreover $x_j \vert [t_0, t_2]$ is injective.
and we have that
$(x_j\vert (t_0, t_2))^{-1} : x_j\left((t_0, t_2)\right) \rightarrow (t_0, t_2)$
is analytic and also its graph is
a semianalytic subset of
$\Bbb{R}\times \Bbb{R}$
and then,
because of lemma \ref{lemaa0}, $(b),$ that there is an extension
$(x_j\vert [t_0, t_2])^{-1} : x_j\left([t_0, t_2]\right) \rightarrow [t_0, t_2]$
which is continuous, and therefore
$x_j \vert [t_0, t_2] : [t_0, t_2] \rightarrow x_j\left([t_0, t_1]\right)$
is an homeomorphism.
Let
$t(s),\, s \in [x_j(t_0), x_j(t_2)],$
be the map
$(x \vert [t_0, t_2])^{-1}_j,$
in other words,
the parameter
$s$
represents the coordinate
$x_j.$
Then we have that
$x([t_0, t_2]) = \{(x_1,...,x_n) : x_i = x_i\left(t(s)\right), s \in
\left[x_j(t_0), x_j(t_1)\right], i = 1,...,n
\}$
and also
$x\left((t_0, t_2)\right) = \{(x_1,...,x_n) : x_i = x_i\left(t(s)\right), s \in
\left(x_j(t_0), x_j(t_1)\right), i = 1,...,n
\}.$
>From this we can easily deduce the assertion of the lemma for the case
of an as-curve
$x :$
$[t_0, t_1) \rightarrow N.$
The case of an as-curve
$x :
(t_0, t_1] \rightarrow N$
can be proven in an entirely similar way.
\quad $\blacksquare$
\\

\begin{lemma}\label{lemmaas2}
Let $x(t),$ $t\in [t_0, t_1)$ be an as-curve in $N,$ which is not a constant. Then there is a $t_2 \in (t_0, t_1)$
such that $x\left([t_0, t_2)\right) - \{x(t_0)\}$ is
nonempty and locally connected at $x(t_0),$
more precisely,
$x\left((t_0, t_3)\right)$
is a neighborhood of
$x(t_0)$
in
$x\left([t_0, t_2)\right) - \{x(t_0)\},$
for all
$t_3 \in (t_0, t_2).$
Moreover,
$x : [t_0,
t_2 ] \rightarrow x \left([t_0, t_2]\right)$
is an homeomorphism, $x \left( (t_0, t_2)\right)$
is an analytic submanifold which is a semianalytic
subset of $N$ and $x : (t_0, t_2) \rightarrow x
\left((t_0, t_2 )\right)$ is an analytic
diffeomorphism.
A similar result holds for an as-curve $x(t),$ $t\in (t_0, t_1],$
in $M.$
\end{lemma}
\textbf{Proof.}
We can show using lemma \ref{lemmaas1}
that there exists $t_2 \in (t_0, t_1)$ such that $x \vert [t_0, t_2]$
is injective, and that $x \vert [t_0, t_2]$
is an homeomorphism onto its image and that
$x : (t_0, t_2 ) \rightarrow x
\left((t_0, t_2 )\right)$ is an analytic diffeomorphism
where
$x\left ((t_0, t_2 )\right)$
is an analytic submanifold.
In particular, we have that
$x(t) \neq x(t_0)$ for all $t \in (t_0, t_2].$
Let
$r > 0$
small be given.
Working in local analytic coordinates centered at
$x(t_0)$ we can show that
continuity of $x(t)$ implies that there exists
$t_r \in (t_0, t_2]$ such that $x\left([t_0,
t_r]\right) \subseteq B_r\left(x(t_0)\right).$ It can be easily
shown that $x\left([t_0, t_r)\right) - \{x(t_0)\} = x\left((t_0,
t_r)\right)$ is connected.
Moreover, for each $s \in [t_0, t_r)$
there is an open ball $B_{\delta}\left(x(s)\right) \subseteq
B_r\left(x(t_0)\right),$ with $\delta = \delta(s),$ such that
$B_{\delta}\left(x(s)\right) \cap x\left([t_r, t_2]\right) = \emptyset$. Then the open set
$$
W = \bigcup_{s \in [t_0, t_r)}B_{\delta}\left(x(s)\right)
$$
satisfies $W \subseteq B_r\left(x(t_0)\right)$ and $W \cap\left(
x\left([t_0, t_2)\right) - \{x(t_0)\}\right)= x\left((t_0,
t_r)\right).$ This shows that $x\left([t_0, t_2)\right)
- \{x(t_0)\}$ is nonempty and locally connected at $x(t_0)$
and also that
$x\left((t_0, t_r)\right)$
is a neighborhood of
$x(t_0)$
in
$x\left([t_0, t_2)\right)
- \{x(t_0)\}.$
Now for each
$t_3 \in (t_0, t_r]$
take
$$
W_{t_3} = \bigcup_{s \in [t_0, t_3)}B_{\bar{\delta}}\left(x(s)\right)
$$
where
$\bar{\delta} = \bar{\delta}(t_3, s)$
and
$B_{\bar{\delta}}\left(x(s)\right)$
satisfies
$B_{\bar{\delta}}\left(x(s)\right)
\subseteq B_r\left(x(0)\right)$
and
$x\left([t_3, t_2]\right) \cap  B_{\bar{\delta}}\left(x(s)\right) = \emptyset.$
Then
$x\left((t_0, t_3)\right)
=
W_{t_3} \cap \left(x\left([t_0, t_2)\right) - \{x(t_0)\}\right).$
This shows that
$x\left((t_0, t_3)\right)$
is a neighborhood of
$x(t_0)$
in
$x\left([t_0, t_2)\right)
- \{x(t_0)\}.$
Then the first case of the lemma is proved.
The case
of an as-curve $x(t),$ $t\in (t_0, t_1],$ in $M$ can be proven in
an entirely similar way.
\quad $\blacksquare$\\

\noindent Inspired by lemma 6.3 of ~\cite{bier88:anal} we will prove the following result about the image and reparametrization of an as-curve.
\begin{lemma}\label{thmas3}
Let $x(t) \in M,$ $t\in [t_0, t_1),$ be an as-curve in $M.$ Then
there is an as-curve $z(s),$ $s \in (s_0 - \delta_1, s_0 + \delta_2),$ in $M,$
for some $\delta_1, \delta_2 > 0$
such that $z(s_0) =
x(t_0),$ $z([s_0, s_0 + \delta_2)) = x ([t_0, t_2))$ for some $t_2 \in (0,
t_1)$ and we also have that $t_2$ and $s_0 + \delta_2$ are such that
$x\vert
[t_0, t_2)$ and $z\vert [s_0, s_0 + \delta_2)$ are homeomorphisms onto
$z([s_0, s_0 + \delta_2)) = x ([t_0, t_2)),$
$x\vert
(t_0, t_2)$ and $z\vert (s_0, s_0 + \delta_2)$ are analytic difeomorphisms onto
$z((s_0, s_0 + \delta_2)) = x ((t_0, t_2))$
which is an analytic submanifold which is also a semianalytic subset and, moreover,
$(x\vert (t_0, t_2))^{-1} \circ (z\vert
(s_0, s_0 + \delta_2)) : (s_0, s_0 + \delta_2) \rightarrow
\Bbb{R}$ is an as-curve in $\Bbb{R}$ which is an analytic
diffeomorfism onto its image $(t_0, t_2).$ Moreover,
$t_2$
can be chosen such that for each
$t_3 \in (t_0, t_2],$
$x\left((t_0, t_3)\right)$ is a
neighborhood of $x(t_0)$ in $x\left([t_0, t_2)\right) - \{x(t_0)\}.$
Similar results hold for as-curves
$x(t)$ in $M,$ where $t\in (t_0, t_1].$
\end{lemma}
\textbf{Proof.}
Using lemma \ref{lemmaas2} and also
lemma 6.3 of
~\cite{bier88:anal}, we can conclude that
there is an as-curve
$z(s),$ $s \in (s_0 - \delta_1, s_0 + \delta_1),$ in $M,$ for some
$\delta_1 , \delta_2 > 0$
such that
$z(s_0) = x(t_0),$
$z ([s_0, s_0 + \delta_2)) =
x ([t_0, t_2))$
for some
$t_2 \in (0, t_1)$
satisfying all the conditions stated in lemma \ref{lemmaas2}.
Since we can also apply lemma \ref{lemaasa}
and
lemma \ref{lemmaas2}
to
$z(s)$
we can also deduce that
$t_2$
and
$\delta_2$
can be chosen  such that
$x\vert
[t_0, t_2)$ and $z\vert [s_0, s_0 + \delta_2)$ are homeomorphisms onto
$z([s_0, s_0 + \delta_2)) = x ([t_0, t_2)),$
$x\vert (t_0, t_2)$
and
$z\vert (s_0, s_0 + \delta_2)$
are analytic diffeomorphisms onto
$z((s_0, s_0 + \delta_2)) = x ((t_0, t_2))$
which is an analytic submanifold which is also a semianalytic subset and, moreover,
$(z\vert (t_0, t_3))^{-1} \circ (x\vert (t_0, t_2)):
(t_0, t_2) \rightarrow \Bbb{R}$
is an as-curve in
$\Bbb{R}$
which is an analytic diffeomorfism onto its image $(s_0, s_0 + \delta_2).$
The rest of the proof can be performed in a similar way.\\
\quad $\blacksquare$
\\

\paragraph{The notion of a pas-curve and of a lcs-curve.}
By definition, a continuous curve in a manifold $M,$ $x(t),$ $t
\in [a, b)$ such that there is a partition $a = t_0 \leq,...,\leq t_m
= b$ such that the restrictions $x \vert [t_i, t_{i + 1}],$ $i =
0,....,m-2,$ $x \vert [t_{m-1}, t_m)$  are as-curves in $M$ is called a
\textsl{continuous piecewise-as-curve in $M$}
(or \textsl{continuous pas-curve in M}).
Each restriction
$x \vert [t_i, t_{i+1}),$
$x \vert (t_i, t_{i+1}],$
$x \vert (t_i, t_{i+1}),$
$x \vert [t_i, t_{i+1}],$
$i = 0,....,m-2,$
$x\vert[t_{m-1}, b),$
$x\vert(t_{m-1}, b)$
is called an as-piece of
$x.$
We define the notion of a \textsl{continuous piecewise-as-curve}
(or \textsl{continuous pas-curve})
$x(t),$ $\,t \in
(a, b],\,t \in [a, b]$ or $t \in (a, b),$ in $M,$
in a similar way.\\

\noindent
We shall also
introduce the notion of a (not necessarily continuous)
\textsl{piecewise-as-curve in $M$}
(or \textsl{pas-curve in $M$})
by eliminating from the
previous definition the condition of continuity at the points
$x(t_i),$
$i = 1,...,m-1,$
belonging to two consecutive as-pieces, and replacing it by the weaker condition of left or right continuity.
More precisely, a pas-curve
$x: [a, b) \rightarrow M$
in
$M$
is defined by the condition
that there is a partition
$a = t_0 \leq,...,\leq t_m = b$
such that, for each
$i = 0,...,m-2$
such that
$t_i \neq t_{i+1},$
the restriction
$x\vert (t_i, t_{i+1})$
is an as-curve in $M$
having compact graph
and
$x\vert (t_{m - 1}, b)$
is an as-curve in $M$
having a continuous extension
$x\vert [t_{m - 1}^+, b)$
obtained by taking the limit at
$t_{m-1}$
on the right. Besides, we require that
$x$
be left continuous or right continuous at each
$t_i,$
$i = 1,...,m-1.$
We define the notion of a \textsl{piecewise-as-curve}
(or \textsl{pas-curve})
$x(t),$ $\,t \in
(a, b],\,t \in [a, b]$ or $t \in (a, b),$ in $M,$
in a similar way.\\

\noindent
We see immediately from this definition that any as-curve is a pas-curve and, moreover, that any
continuous pas-curve is a pas-curve, in a natural way.
Using lemma \ref{lemmasuban1} we can conclude that, given a pas-curve
$x : [a,b)\rightarrow M$
as before, for each $i = 0,...,m-2$
there is a well defined as-curve denoted
$x\vert [t_i^+, t_{i+1}^-],$
which, for
each
$i = 0,...,m-2$
such that
$t_i \neq t_{i+1},$
is the uniquely defined continuous extension of the restriction
$x\vert (t_i, t_{i+1}).$
We can easily see from lemma \ref{lemmasuban1} that
if the graph of
$x\vert [t_{m - 1}^+, b)$
is not relatively compact then it must be closed.
Similar statements hold for pas-curves
$x(t),$ $\,t \in
(a, b],\,t \in [a, b]$ or $t \in (a, b),$ in $M.$
\\

\noindent
Given a pas-curve
$x(t),$
$t \in [a, b),$
as above, we have that some of the restrictions
$x\vert [t_i, t_{i + 1}),$
$x\vert (t_i, t_{i + 1}],$
$x\vert [t_i, t_{i + 1}],$
$x\vert (t_i, t_{i + 1}),$
for each
$i = 0,...,m-2,$
and some of the restrictions
$x\vert [t_{m-1}, b)$
$\vert (t_{m-1}, b),$
satisfy the required continuity condition at the closed end of the interval and then they
are as-curves.
Each one of those restrictions which is an as-curve is called an
\textsl{as-piece} of
$x.$
If the graph of
$x$
is relatively compact
we will call each as-curve
$x\vert [t_i^+, t_{i+1}^-],$
$i = 0,...,m-2,$
$x\vert [t_{m-1}^+, b^-]$
an \textsl{as-piece closure}
of
$x.$
If
$\operatorname{graph} x$
is not relatively compact then the \textsl{as-piece closures} are
$x\vert [t_i^+, t_{i+1}^-],$
$i = 0,...,m-2,$
$x\vert [t_{m-1}^+, b).$
We must observe that an as-piece closure is an as-curve but it is not always an as-piece
of
$x.$
By definition, the \textsl{open as-pieces}
are
$x\vert (t_i, t_{i+1}),$
$i = 0,...,m-2,$
$x\vert (t_{m-1}, b).$
Similar definitions hold for curves
$x(t),$
where
$t \in (a, b],$ $t \in [a, b],$ $t \in (a, b).$\\

\noindent
We need the following definition.
A pas-curve
$y$
in
$M$
is a \textsl{refinement}  of a pas-curve
$z$
in $M$
if
$\operatorname{graph}y \ = \operatorname{graph}z$
and
the graph of
each as-piece of
$y$
is contained in the graph of some as-piece of
$z.$\\

\noindent We introduce the notation
$C^{as}([a, b), M)$
to denote the set of all
as-curves
$x : [a, b) \rightarrow M,$
where
$M$
is a given manifold.
In a similar way, we define
$C^{as}((a, b], M),$
$C^{as}((a, b), M),$
$C^{as}([a, b], M),$
and also
$C^{pas}([a, b), M),$
$C^{pas}((a, b], M),$
$C^{pas}((a, b), M),$
$C^{pas}([a, b], M).$

\begin{lemma}\label{lemmaab}
$(a)$ Every pas-curve in the manifold $M,$ say $x(t),$ $t
\in [a, b),$ ($t \in (a, b],$ $t \in [a, b],$ $t \in (a, b)$), is
a subanalytic curve in $M.$\\
\noindent
$(b)$ Let $x(t),$ $t \in [a, b),$ be a subanalytic curve
in $M$ whose graph is a relatively compact subset of
$\Bbb{R}\times M$
and which is left or right continuous at each
$t \in [a, b).$
Then there is a, not necessarily unique,
pas-decomposition of $x(t),$
that is, a partition of $[a, b],$ say $a = t_0 \leq t_1,...,\leq t_r = b,$
such that $x\vert [t_i^+, t_{i + 1}^-],$ $i = 0,...,r-2$ and $x\vert
[t_{r-1}^+, b)$ are as-curves. A similar result holds for
subanalytic curves $x(t),$ $t \in (a, b],$ $t \in [a, b]$  or $t
\in (a, b)$.
\end{lemma}
\textbf{Proof.} The proof of $(a)$ is a consequence of the fact that
each as-piece of a pas-curve has a semianalytic graph and also the
fact that a finite union of semianalytic subsets is a semianalytic subset.
To prove $(b)$ we observe first that, because of lemma
\ref{lemmasuban1}, there is a continuous extension,
$\bar{x}(t),$
$t \in [a, b],$
of
$x(t).$
Also lemma \ref{lemmasuban1} implies that for each
$t \in [a, b]$
the right and left limits
$x(t^-),$ $x(t^+)$
exists.
Then, by lemma \ref{lemaasa},
for each
$t\in (a, b)$
there exists
$\epsilon > 0$
such that
$x\vert [t^+, (t + \epsilon)^-],$
$x\vert [(t + \epsilon)^+, t^-]$
are
as-curves,
and also
$x\vert [a^+, (a + \epsilon)^-],$
$x\vert [(b - \epsilon)^+, b^-]$
are as-curves for some
$\epsilon > 0.$
Then the proof follows by a
standard compactness argument.
\quad $\blacksquare$\\

\noindent Since we want to work within the subanalytic category,
it seems that a good choice for a class of possible solutions to a
given IDE would be the class of subanalytic curves in $M$ whose
graph is a relatively compact subset and which are left or right
continuous at each point. The decomposition, predicted
in lemma \ref{lemmaab}, of such a curve in as-pieces, which gives
its structure as a pas-curve in $M,$ which is not unique, is what
will enable us to define, in the next paragraph, how such curves
can be interpreted as being pas-solutions to a given IDE. In other
words, one convenient class of possible solutions to a given IDE
would be the class of subanalytic curves having a relatively
compact graph, interpreted as pas-curves, according to lemma
\ref{lemmaab}.
As another choice for a class of possible solutions one could choose the class of all pas-curves in $M.$
We are going to work also with a bigger class of
curves, namely the class of curves that are {\sl
locally-compact-subanalytic}, or {\sl lcs-curves}, in the sense that their restriction
to any compact subinterval is a subanalytic curve whose graph is
compact
and which are left or right
continuous at each point.
More precisely, we define
\begin{definition}\label{definlcs}
Let $M$ be a given manifold. Then we define the following classes of
locally-compact-subanalytic curves, also called lcs-curves:
\\
\noindent $(a)$ $C^{lcs}\left([t_0, t_1), M\right),$ where $t_1 \in
\Bbb{R}$ or $t_1 = +\infty$ is the set of all functions $x : [t_0,
t_1) \rightarrow M$ such that $x\vert [t_0, t_2]$ is a pas-curve,
for each $t_2 \in (t_0, t_1).$\\
\noindent
$(b)$ $C^{lcs}\left((t_0, t_1], M\right),$ where $t_0 \in
\Bbb{R}$ or $t_0 = -\infty$ is the set of all functions $x : (t_0,
t_1] \rightarrow M$ such that $x\vert [t_2, t_1]$ is a pas-curve,
for each $t_2 \in (t_0, t_1).$
\\
\noindent
$(c)$ $C^{lcs}\left((t_0, t_1), M\right),$ where $t_0 \in
\Bbb{R}$ or $t_0 = -\infty$ and $t_1 \in \Bbb{R}$ or $t_1 =
+\infty$ is the set of all functions $x : (t_0, t_1) \rightarrow M$
such that $x\vert [\bar{t}_0, \bar{t}_2]$ is a pas-curve, for each
$[\bar{t}_0, \bar{t}_1] \subseteq (t_0, t_1).$\\
\noindent
$(d)$ $C^{lcs}\left([t_0, t_1], M\right),$ where $t_0 \in
\Bbb{R}$ and $t_1 \in \Bbb{R}$ is the set of all functions $x : [t_0, t_1] \rightarrow M$
such that $x$ is a pas-curve.
\end{definition}
We may define the notion of an \textsl{as-piece} and also the notion of a \textsl{refinement}
of a given
lcs-curve in a similar way as we did in the case of a pas-curve.\\

\noindent
{\bf Remark.}
(i) With the notation of definition \ref{definlcs}, we have, respectively for
$(a),$ $(b),$ $(c)$ and $(d),$
that each one of the graphs, $\operatorname{graph}(x\vert [t_0,
t_2]),$ $\operatorname{graph}(x\vert [t_2, t_1]),$
$\operatorname{graph}(x\vert [\bar{t}_0, \bar{t}_1])$
and
$\operatorname{graph}(x\vert [t_0, t_1])$
is compact.\\
\noindent (ii) In
view of lemma \ref{lemmaab} we may replace {\sl pas-curve} by
{\sl subanalytic curve
which is left or right continuous at each point}
in definition \ref{definlcs} and we will obtain
an equivalent definition. The decomposition of each subanalytic curve $x\vert [t_0,
t_2],$ $x\vert [t_2, t_1],$ $x\vert [\bar{t}_0, \bar{t}_1]$
or
$x\vert [t_0, t_1]$
as a
pas-curve in $M$ is not unique.
For a given lcs-curve in $M$ there is a decomposition in at most a countable number of
as-pieces, which is not unique.\\

\paragraph{Solutions to IDE: as-solutions, pas-solutions and lcs-solutions.}
Now we introduce the notion of solution to a given IDE which is
convenient for the purposes of the present paper.
\begin{definition}
$(a)$ An \textsl{as-solution} $x(t),$ $t\in [t_0, t_1)$ in $M$ to a
given IDE $(a , f)$ in $M$ is an as-curve in $M$ which satisfies
$(a, f)$ for all $t \in (t_0, t_1),$ that is,
$a\left(x(t)\right)\dot{x}(t) = f\left(x(t)\right),$ for all $t
\in (t_0, t_1).$ Similar statements hold for as-solutions $x(t)$
in $M$ where $t\in (t_0, t_1],$
$t \in (t_0, t_1)$
or
$t \in [t_0, t_1]$
with
$t_0 \neq t_1.$
\\
\noindent
$(b)$ A \textsl{pas-solution} in $M$ to a given IDE is a pas-curve in $M$
such that each open nonempty as-piece
is an as-solution in $M.$\\
\noindent
$(c)$ A \textsl{lcs-solution} in $M$ to a given IDE is a curve in $M$
belonging to
$C^{lcs}\left([t_0, t_1), M\right),$
$C^{lcs}\left((t_0, t_1], M\right),$
$C^{lcs}\left((t_0, t_1), M\right)$
or
$C^{lcs}\left([t_0, t_1], M\right)$
such that each restriction to a compact subinterval has a pas-decomposition which is a pas-solution in
$M.$
\end{definition}
\paragraph{Lifted and projected solutions.}
We have the following result.
\begin{lemma}\label{lemmaas3}
$(a)$ Let $y(t),$ $t\in
[t_0, t_1),$ ($t\in (t_0, t_1],$ )
be a given as-solution
in
$M^k$
to the system $(a_k , f_k)$
described in the previous section,
for some
$k = 1, 2,...q.$
Then for each  $t_2 \in (t_0,
t_1),$ $y(t)$ is projected via $\pi_{k-1}$ into an as-solution
$x(t)$
to the system
$(a_{k-1}, f_{k-1}),$
$x(t) = \pi_{k-1} \left(y(t)\right),$ $t\in[t_0, t_2],$
(respectively, $t\in [t_2, t_1]$ ),
in $M^{k-1}.$
\\
\noindent
$(b)$ Assume that we have an as-solution $y(t),$
$t\in [t_0, t_1),$ ($t\in (t_0, t_1]$ ), in $M^k,$
for some $k = 1, 2,...q,$
to the system $(a_k , f_k)$ described in the previous section. Then
for each $s = 0,...,k-1,$ $y(t)$ is projected via
$\pi_s\circ...\circ \pi_{k-1}$ into an as-solution $x(t)$ to the system
$(a_s, f_s),$
$x(t) = \pi_s\circ...\circ \pi_{k-1} \left(y(t)\right),$
$t\in[t_0, t_2],$
(respectively, $t\in [t_2, t_1],$ ), in $M^s,$ for each $t_2 \in (t_0, t_1).$
\end{lemma}
\textbf{Proof.} Part $(a)$ is easy to prove using the fact that $\pi_{k-1}$
is an analytic map,
and also lemma \ref{lemmaas0}, $(c)$.
Part $(b)$ follows using $(a).$
\quad $\blacksquare$
\\

\noindent
>From the previous lemma we can deduce that if
$y(t),$
$t \in [t_0, t_1)$
is an as-solution to
$(a_k, f_k)$
in
$M^k,$
for some
$k = 1,...,q$
and
$x(t) = \pi_s \circ,...,\pi_k \left(y(t)\right)$
is an as-curve in
$M^s,$
for some
$s = 1,...,k,$
then
$x(t),$
$t \in [t_0, t_1)$
is an as-solution to
$(a_s, f_s).$
We will call $y(t)$ a
\textsl{lifted as-solution of}
$x(t),$
and
$x(t)$ the
\textsl{projected as-solution of} $y(t).$
A similar definition holds for the case of as-solutions
$y(t),$
$t \in (t_0, t_1],$
$t \in (t_0, t_1),$
$t \in [t_0, t_1].$
The notions of a {\sl
projected} and {\sl lifted pas-solution}
are defined as follows.
A pas-solution to $(a_k, f_k),$ say $y(t)\in M^k,$
$t \in [t_0, t_1)$
is a \textsl{lifted pas-solution  of}
a pas-solution $x(t)$ to $(a_s, f_s)$
if each as-piece of
$y(t)$ is projected, via the composition
$\pi_s \circ...\circ\pi_{k-1},$
onto an as-piece of $x(t) = \pi_s \circ...\circ\pi_{k-1} y(t),$
which is then an as-solution
to
$(a_s, f_s).$
We call
$x(t)$
the \textsl{projected pas-solution
of}
$y(t).$
A similar definition holds for the case of pas-solutions
$y(t),$
$t \in (t_0, t_1],$
$t \in (t_0, t_1),$
$t \in [t_0, t_1].$
\paragraph{Main results.}

\begin{theorem}\label{liftheorem}
$(a)$ Let $y(t),$  $t \in [t_0, t_1)$ (respectively, $t \in (t_0,
t_1],$) be an as-solution to $(a_k, f_k)$ in $M^{k},$ $k =
1,...,q.$ Then $x(t) = \pi_{k-1} \left(y(t)\right),$ $t \in [t_0,
t_2]$ (respectively, $t \in [t_2, t_1]$) is an as-solution to
$(a_{k-1}, f_{k-1})$ in $M^{k-1},$ for each $t_2 \in (t_0,
t_1).$
The previous statement holds true
if we replace
an \textsl{as-solution} by a \textsl{pas-solution} or a \textsl{lcs-solution}.
\\
\noindent
 $(b)$ If $x(t),$ $t \in [t_0, t_1)$ (respectively, $t \in
(t_0, t_1]$) is an as-solution to $(a_{k-1}, f_{k-1})$ in
$M^{k-1}$
such that
$x(t) \in M^{k-1}_0,$ $t \in [t_0, t_1)$
(respectively, $t \in
(t_0, t_1]$),
$k = 1,...,q$ then there exists $t_2\in (t_0, t_1)$ and
a lifted as-solution $y(t),$ $t \in [t_0, t_2]$
(respectively, $t \in
[t_2, t_1]$)
of
$x\vert [t_0, t_2]$
(respectively, $x\vert [t_2, t_1]$)
to $(a_k, f_k)$ in $M^k,$ in particular, $x(t) =
\pi_{k-1}\left(y(t)\right),$
$t\in [t_0, t_2]$ (respectively, $t\in [t_2, t_1].$)
\\
\noindent
$(c)$ Let $x(t),$ $t \in [t_0, t_1]$ be an as-solution to
$(a_{k-1}, f_{k-1})$ in $M^{k-1}$
such that
$x(t) \in M^{k-1}_0,$ $t \in [t_0, t_1],$
$k = 1,...,q.$ Then there is a
lifted pas-solution of a refinement of
$x,$ say $y(t)$, $t \in [t_0, t_1]$ to $(a_k, f_k)$ in
$M^k,$ in particular, $x(t) = \pi_{k-1}\left(y(t)\right),$
$t\in [t_0, t_1].$
Moreover, if $C
\subseteq \pi_{k-1}^{-1}(x([t_0, t_1])),$ $k = 1,...,q,$
is a compact semianalytic subset
such that
$\pi_{k-1} (C) = x([t_0, t_1])$ then there is a lifted
pas-solution
of a refinement of $x,$ say
$y(t)\in C,$ $t \in [t_0, t_1],$ to $(a_k, f_k)$ in
$M^k,$ in particular, $x(t) = \pi_{k-1}\left(y(t)\right),$ $t \in [t_0, t_1].$
\\
\noindent
$(d)$ Let $x(t),$ $t \in [t_0, t_1]$ be a pas-solution to
$(a_{k-1}, f_{k-1})$ in $M^{k-1}$
such that
$x(t) \in M^{k-1}_0,$ $t \in [t_0, t_1],$
$k = 1,...,q.$ Then there is a
lifted pas-solution of a refinement of
$x,$ say  $y(t),$ $t \in [t_0, t_1]$ to $(a_k, f_k)$ in
$M^k,$ in particular, $x(t) = \pi_{k-1}\left(y(t)\right),$
$t\in [t_0, t_1].$
Moreover, if $C \subseteq \pi_{k-1}^{-1}(x([t_0,
t_1])),$ $k = 1,...,q,$
is a compact semianalytic subset
such that $\pi_{k-1} (C) = x([t_0, t_1]),$
then there is a lifted pas-solution
of a refinement of
$x,$
say
$y(t)\in C,$ $t \in
[t_0, t_1],$ to $(a_k, f_k)$ in $M^k,$ in particular, $x(t) =
\pi_{k-1}\left(y(t)\right),$
$t\in [t_0, t_1].$
\\
\noindent
$(e)$ Let $x(t),$ $t \in [t_0, t_1),$ ($t \in (t_0, t_1],$
$t \in [t_0, t_1],$ $t \in (t_0, t_1)$) be a lcs-solution to $(a_{k-1}, f_{k-1})$
in $M^{k-1},$
such that
$x(t) \in M^{k-1}_0,$ $t \in [t_0, t_1),$
(respectively, $t \in (t_0, t_1],$
$t \in [t_0, t_1],$ $t \in (t_0, t_1)$),
$k = 1,...,q.$
Then there is a lifted lcs-solution
of a refinement of
$x,$ say
$y(t),\,
t \in [t_0, t_1),$ (respectively $t \in (t_0, t_1],$ $t \in [t_0, t_1]$
or
$t \in (t_0, t_1)$) to $(a_k, f_k)$ in $M^k,$ in particular, $x(t) =
\pi_{k-1}\left(y(t)\right),$
$t \in [t_0, t_1),$ ($t \in (t_0, t_1],$
$t \in [t_0, t_1],$ $t \in (t_0, t_1)$).
Moreover, if $C \subseteq
\pi_{k-1}^{-1}(x([t_0, t_1))),$ (respectively, $C \subseteq
\pi_{k-1}^{-1}(x((t_0, t_1])),$
$C \subseteq
\pi_{k-1}^{-1}(x([t_0, t_1])),$
$C \subseteq
\pi_{k-1}^{-1}(x((t_0, t_1)))$), $k = 1,...,q,$ is a subset such that
$\pi_{k-1} (C) = x([t_0, t_1)),$ (respectively, $\pi_{k-1}(C) =
x((t_0, t_1]),$ $\pi_{k-1}(C) =
x([t_0, t_1]),$ $\pi_{k-1}(C) = x((t_0, t_1))$), and, besides,
$C \cap \pi_{k-1}^{-1}(x([\bar{t}_0, \bar{t}_1]))$ is compact subanalytic
for each compact subinterval $[\bar{t}_0, \bar{t}_1] \subseteq
[t_0, t_1)$ (respectively, $[\bar{t}_0, \bar{t}_1] \subseteq (t_0,
t_1],$
$[\bar{t}_0, \bar{t}_1] \subseteq [t_0,
t_1],$
$[\bar{t}_0, \bar{t}_1] \subseteq (t_0, t_1)$), then
there is a lifted lcs-solution
of a refinement of
$x,$ say
$y(t)\in C,$ $t \in
[t_0, t_1),$ (respectively $t \in (t_0, t_1],$
$t \in [t_0, t_1],$
$t \in (t_0,
t_1)$) to $(a_k, f_k)$ in $M^k,$ in particular, $x(t) =
\pi_{k-1}\left(y(t)\right),$
$t \in [t_0, t_1),$ ($t \in (t_0, t_1],$
$t \in [t_0, t_1],$ $t \in (t_0, t_1)$).
\end{theorem}
\textbf{Proof.}
We are going to give a detailed proof of the case
$k = 1$
only, since the cases
$k = 2,...,q$
can be proven in an entirely similar way.
Part $(a)$ is an easy consequence of lemma \ref{lemmaas3}.
In order to prove $(b)$, $(c)$, $(d)$ and $(e)$ we are going to prove first several facts,
namely,
$(i),$ $(ii),$ $(iii),$ $(iv)$ and $(v)$ below. These facts will be proven under the assumption that $x(t)\in M_0,$ $t\in [t_0, t_1],$ is
an as-solution to $(a, f)$ in $M,$ the curve $x(t)$ is
simple, that is, $x(a) \neq x(b)$ for all $a, b \in [t_0, t_1]$
such that $a \neq b,$ and moreover,
$x: [t_0, t_1] \rightarrow x([t_0, t_1])$ is
an homeomorphism
and also
$x((t_0, t_1))$
is an analytic submanifold and
$x\vert (t_0, t_1) : (t_0, t_1) \rightarrow x((t_0, t_1))$
is an analytic diffeomorphism.
We can assume without loss of generality (for instance, using
Whitney embedding theorem) that $M\subseteq U$ and $M^1 \subseteq
V$ are analytic submanifolds of $U$ and $V,$ where $U$ and $V$
are real finite dimensional vector spaces.\\

\noindent
 $(i)$ The map $\pi_0$ can be described as the
restriction $p\vert (\operatorname{graph}\pi_0)$ to
$\operatorname{graph}\pi_0 \subseteq V\times U$ of the projection
onto the second factor $p : V\times U \rightarrow U.$
Since
$p\vert (\operatorname{graph}\pi_0)$ is a proper analytic map we
have that $\left(p\vert (\operatorname{graph}\pi_0)\right)^{-1}
\left(x\left([t_0, t_1]\right)\right)$ is a compact semianalytic
subset of $V\times U,$ therefore, according to lemma 3.4 of
~\cite{bier88:anal} it is a finite union of
relatively compact
connected smooth
semianalytic subsets $A$ such that for each $A$
$\operatorname{rank}\left( p\vert A\right)$ is constant on $A.$
Since $\operatorname{dim}\left(x([t_0, t_1]\right) = 1,$ it is
easy to see that for each $A$ $\operatorname{rank}\left( p\vert
A\right)$ is 0 or 1, and moreover, $x([t_0, t_1])$ is the union of
those $p(\bar{A})$ such that $p\vert A$ has rank 1 and therefore
$p(\bar{A})$ is not a point.
We observe that if there is a
compact semianaliytic subset
$C \subseteq \pi_0^{-1}
\left(x([t_0, t_1])\right)$
such that
$\pi_0 (C) = x([t_0, t_1])$
then we have that
$C_1 = (p\vert \operatorname{graph} \pi_0)^{-1} \left(x([t_0, t_1])\right) \cap (C\times U)$
is compact and semianalytic.
According to
lemma 3.4 of ~\cite{bier88:anal}, $C_1$ is a finite union of
relatively compact
connected smooth semianalytic subsets $A$ having the same
properties as before.
\\

\noindent $(ii)$ Let $\alpha \in [t_0, t_1]$ be fixed. Then there
is an $A,$ say $A = A_0,$ such that $x(\alpha) \in p(\bar{A}_0).$
Since $\bar{A}_0$ is connected and compact
and the curve $x(t),$\, $t \in [t_0, t_1],$ is simple, we have
that $p(\bar{A}_0)$ is homeomorphic to a closed interval (possibly
of zero length) via the curve $x,$ say $p(\bar{A}_0) =
x\left([a_1, a_2]\right),$ where $[a_1, a_2]\subseteq [t_0, t_1].$
We can assume without loss of generality that $p(\bar{A}_0)$ is
homeomorphic to a closed interval of nonzero length.
\\

\noindent
 $(iii)$ By using theorem 6.10 of ~\cite{bier88:anal} we
can see that for given points $q_i \in
p^{-1}\left(x(a_i)\right)\cap \bar{A}_0,$ $i = 1, 2,$ so in
particular, we have $q_1 \neq q_2,$ there is a continuous
semianalytic curve $w(s),$ $s \in [s_0, s_1],$ in $\bar{A}_0$ such
that $w(s_i) = q_{i+ 1},$ $i = 0, 1.$ Then $w([s_0, s_1])$ is a
compact semianalytic subset of $V\times U$ of dimension 1 and we
have $p\left(w([s_0, s_1])\right) = x([a_1, a_2])].$ Using lemma
3.4 of ~\cite{bier88:anal} we see that since $w([s_0, s_1])$ has
dimension 1 it is a finite union of relatively compact
connected smooth semianalytic
subsets $B$ of dimension less or equal than 1, such that the
restriction of the projection $p\vert B$ has constant rank of
value 0 or 1. Since each $\bar{B}$ is connected and compact we
have that $p\left(\bar{B}\right)$ is homeomorphic to a closed
interval, say $p\left(\bar{B}\right) = x\left([a_B, b_B]\right).$
For at least some $B$ such that $p\left(\bar{B}\right)$ is not a
point one must have that $p\left(\bar{B}\right)$ contains the
point $x(\alpha).$
Observe that
$p\left(B\right) = x\left((a_B, b_B)\right)$
is an analytic submanifold which is a subanalytic
subset and that
$p\vert B : B \rightarrow  x\left((a_B,
b_B)\right)$
is an analytic diffeomorphism. This gives, in
particular, a parametrization of the analytic submanifold $B$ with the
parameter $t,$ namely, $z(t) = \left(p\vert B\right)^{-1}
\left(x(t)\right),$ $t \in (a_B, b_B).$
We can give a definition of
$\operatorname{graph}\left(\left(p\vert B\right)^{-1} \circ x \right)$
as a subanalytic subset of
$\Bbb{R}\times M^1 \times M$
as follows.
We have
$\operatorname{graph}\left(\left(p\vert B\right)^{-1} \circ x \right)
=
\{(t, z) \in \Bbb{R}\times M^1 \times M: z \in B ,
p(z) = x, (t, x) \in \operatorname{graph} (x \vert (a_B, b_B))\},$
which defines
$\operatorname{graph} \left(\left(p\vert B\right)^{-1} \circ x \right)$
by subanalytic conditions, since
$B$ is a semianalytic subset of
$M^1 \times M,$
$\operatorname{graph} x$
is a subanalytic subset of
$\Bbb{R}\times M,$
by definition, and
$\operatorname{graph}(x \vert (a_B, b_B))
=
\{(t, x) \in \Bbb{R}\times M:
(t, x) \in \operatorname{graph} x, t\in  (a_B, b_B)\}.$
Then using lemma \ref{lemmasuban1}
we can deduce that there is a uniquely defined continuous extension,
which we will call
$z$
by a slight abuse of notation,
$z(t),$  $t \in [a_B, b_B],$
whose image is the semianalytic subset
$\bar{B} = z([a_B, b_B]),$
which is an as-curve in
$M^1 \times M.$
We have, in particular, that
$z(t),$
$t \in [a_B, b_B]$
is an as-curve such that
$p \left(z(t)\right)
=
x(t),$
$t \in [a_B, b_B].$
It is clear that the extension
$z(t),$
$t \in [a_B, b_B]$
is given by
$\left(p\vert \bar{B}\right)^{-1} \circ x$
\\

\noindent
 $(iv)$ Assume that $x(\alpha) \in p\left(B\right).$
We have the as-curve
$z(t) = \left(y(t), x(t)\right),$
$t \in [a_B, b_B]$
in $\operatorname{graph}(\pi_0),$ therefore $p \left(z(t)\right) =
x(t),$ for $t \in [a_B, b_B],$
then the curve $y(t)$ satisfies $\pi_0
\left(y(t)\right) = x(t),$
$t \in [a_B, b_B].$
Using this it becomes clear from the definition of
$(a_1, f_1)$ that $y(t)$ satisfies the system $(a_1, f_1),$ for $t
\in (a_B, b_B).$
We have that, for any
$[\alpha - \epsilon_1, \alpha + \epsilon_2] \subseteq [a_B, b_B],$
where
$\epsilon_1, \epsilon_2 \geq 0,$
$y(t),$ $t \in [\alpha - \epsilon_1, \alpha + \epsilon_2]$
is an as-curve in $M^1.$
In fact, this is a direct consequence of lemma
\ref{lemmaas0}, $(c)$, since the projection
$q : M^1 \times M \rightarrow M^1$
is an analytic map.
It is clear that $y(t)$ satisfies the system $(a_1, f_1),$ for $t \in
(\alpha - \epsilon_1, \alpha + \epsilon_2).$
\\

\noindent
$(v)$ Assume now that $x(\alpha) \in
p\left(\bar{B}\right) - p\left(B\right),$ then $\alpha = a_B$ or
$\alpha = b_B.$ If $\alpha = a_B$ (respectively $\alpha = b_B$) we
can proceed in a similar way as we did in $(iv)$ and we have an
as-curve $z(t) = \left(y(t), x(t)\right),$ $t \in [\alpha,
\alpha + \epsilon_2]$ (respectively, $t \in [\alpha - \epsilon_1,
\alpha]$) in $\operatorname{graph}(\pi_0),$ then,
in particular,
$p\left(z(t)\right) = x(t),$ $ t \in [\alpha,
\alpha + \epsilon_2],$ (respectively, $t \in [\alpha - \epsilon_1,
\alpha)]$.
Then the curve $y(t) = q\left(z(t)\right)$ satisfies
$\pi_0 \left(y(t)\right) = x(t),$ $t \in [\alpha, \alpha +
\epsilon_2]$ (respectively, $t \in [\alpha - \epsilon_1, \alpha]$)
and is an as-curve.
It is clear that $y(t)$ satisfies the system $(a_1, f_1),$ for $t \in
(\alpha, \alpha + \epsilon_2)$ (respectively $t \in (\alpha -
\epsilon_1, \alpha)$).
\\

\noindent We are going to prove $(b)$. If $x(t)  = x(t_0)$ is a
constant then it can be lifted to a constant curve $y(t) =
y(t_0),$ where $y(t_0) \in \pi_0^{-1}\left(x(t_0)\right),$ which
solves the problem, so let us assume that $x(t)$ is not a
constant. By conveniently lowering the value of $t_1$ we can assume
without loss of generality that $x(t)\in M_0,$ $t\in [t_0, t_1],$
is an as-solution to $(a, f)$ in $M.$
Moreover, by lowering the value of
$t_1$
if necessary and
using lemma
\ref{lemmaas1} we can assume that the curve $x(t)$ is simple, that
is $x(a) \neq x(b)$ for all $a, b \in [t_0, t_1]$ such that $a
\neq b,$ and moreover, that $x : [t_0, t_1] \rightarrow x([t_0, t_1])$ is
an homeomorphism onto
$x([t_0, t_1]),$
$x \left((t_0, t_1)\right)$
is an analytic submanifold of
$M$
and
$x \vert (t_0, t_1) : (t_0, t_1) \rightarrow x\left((t_0, t_1)\right)$
is an analytic diffeomorphism.
By using $(v)$ with $\alpha = t_0,$ we must have
$a_B = t_0$ and then the proof of $(b)$ follows from $(v)$ by taking
$\alpha + \epsilon_2 = t_2.$
The case of an interval
$(t_0, t_1]$
can be proved in an entirely similar way.
\\

\noindent We are going to prove $(c).$
As in the proof of
$(b),$
the case in which
$x$
is constant is trivial, so we shall assume that
$x$
is not a constant.
Using
$(b)$ we can conclude that for each
$\bar{t} \in (t_0, t_1)$
there are intervals
$[\bar{t} - \epsilon_1, \bar{t}]
\subseteq [t_0, t_1]$
and
$[\bar{t}, \bar{t} + \epsilon_2] \subseteq [t_0, t_1],$
where
$\epsilon_1, \epsilon_2 > 0,$
and lifted as-solutions
$y^{-}_{\bar{t}}(t),$
$t \in [\bar{t} - \epsilon_1, \bar{t}],$
and
$y^{+}_{\bar{t}}(t),$
$t \in [\bar{t}, \bar{t} + \epsilon_2],$
in
$M^1,$
of
$x \vert [\bar{t} - \epsilon_1, \bar{t}]$
and
$x \vert [\bar{t}, \bar{t} + \epsilon_2]$
to
$(a_1, f_1).$
Also
if
$\bar{t} = t_0$
($\bar{t} = t_1$)
there is a lifted as-solution
$y^{+}_{t_0}(t)$
of
$x \vert [t_0, t_0 + \epsilon_2]$
(respectively,
$y^{-}_{t_1}(t)$
of
$x \vert [t_1 - \epsilon_1, t_1],$)
to
$(a_1, f_1).$
Using this and a compactness argument we can conclude that there is a partition of
$[t_0, t_1],$
say,
$t_0 = \bar{t}_0 < \bar{t}_1,...,< \bar{t}_r = t_1,$
such that there are lifted as-solutions
$y_i(t),$
$t \in [\bar{t}_i, \bar{t}_{i + 1}],$
of
$x\vert [\bar{t}_i, \bar{t}_{i + 1}],$
to
$(a_1, f_1),$
$i = 0,...,r-1.$
Then a lifted pas-solution can be obtained
by gluing together the lifted solutions
$y_i(t),$
$i = 0,...,r-1$
into a single pas-solution
$y(t), t \in [t_0, t_1],$
which is always possible, and in more than one way if
$r > 0.$
If there is a set
$C$
satisfying the conditions stated in
$(c)$
then we can prove that there is a lifted pas-solution satisfying the required conditions proceeding as above
and taking into account the comments made at the end of $(i).$
\\

\noindent The proof of $(d)$ can be obtained as a consequence of $(c)$
by conveniently gluing together a finite number of lifted pas-solutions each one corresponding to each as-piece
of the given pas-solution.\\

\noindent The proof of $(e)$ can be obtained as a consequence of $(d)$
by conveniently gluing together conveniently chosen at most countable lifted pas-solutions.
\quad $\blacksquare$\\

\noindent
The next theorem
adds some more information and completes the picture
of the relationship between solutions and lifted
solutions to a given IDE.\\

\noindent Let $x(t),$ $t \in [\alpha, \beta]$ be an as-solution to
$(a, f)$ in $M^k,$
for some
$k = 0,...,q-1.$
If for some
$s_0\in [\alpha, \beta]$
we have
$x(s_0) \in M^k_2$
then we must have that there is at most a
finite number of
$t \in [\alpha, \beta],$
say
$\alpha \leq t_1
<....<t_r \leq \beta,$
such that
$x(t_i) \in M^k_0.$
This is because
$M^k_0$
is defined by analytic equations.
In this case we will call
$x$
a
\textsl{normal as-solution to}
$(a_k, f_k)$
\textsl{in}
$M^k.$
If
$x(t),$ $t \in [\alpha, \beta]$ is a pas-solution to
$(a_k, f_k)$ in $M^k,$
for some
$k = 0,...,q-1$
such that all its as-pieces are normal we will call
$x$
a \textsl{normal pas-solution to}
$(a_k, f_k)$
\textsl{in}
$M^k.$\\

\noindent
Let us introduce the following notation.
$\tilde{M} = \bigsqcup_{k = 0}^{q}M^k,$
where
$\bigsqcup$
means disjoint union,
$\tilde{\pi} : \tilde{M} \rightarrow M,$
where
$\tilde{\pi}(x) = \pi_k(x)$
if
$x \in M^k.$
Now assume that we have a finite collection of normal pas-solutions of the
following type. For each
$i = 1,...,r$
let
$y_i(t),$ $t \in
[t_{i,0}, t_{i,1}]$
be a normal pas-solution to some of
the IDE $(a_k, f_k),$ $k = 0,...,q,$
say
$(a_{k_i},
f_{k_i}),$
where
$k_i \in \{0,...,q\},$
$i = 1,...,r,$
so that
$y_i(t) \in M^{k_i}_0$
only for at most a finite number of
$t \in [t_{i,0}, t_{i,1}].$
Assume, without loss of generality, that
$t_{i,1} =
t_{i + 1,0},$
for
$i = 1,...,r-1,$
so we have a partition
$t_{1,0}
< t_{1,1} = t_{2,0} <...<t_{r,1}$
of the interval
$[t_{1,0},
t_{r,1}].$
Then the collection of curves
$y_i,$
$i = 1,..,r$
may be thought of as a single pas-curve
in
$\tilde{M},$
say
$y(t),$
$t \in [t_{1,0},
t_{r,1}].$
We will say that
$y(t),$
$t \in [t_{1,0},
t_{r,1}]$
is a \textsl{normal pas-solution} to
$(\tilde{a}, \tilde{f})$
\textsl{in}
$\tilde{M}.$
It follows from theorem
\ref{liftheorem}, $(a)$
that
$x(t),$ $t \in [t_{1,0},\,t_{r,1}],$ given by $x(t) =
\pi_0 \circ...\circ \pi_{k_i - 1} y_i(t),$ $t \in
[t_{i,0},\,t_{i,1}],$
$i = 1,..,r,$
or, equivalently,
$x(t) = \tilde{\pi}\left(y(t)\right),$
$t \in [t_{1,0},
t_{r,1}],$
is a
pas-solution to
$(a, f).$
In this situation we will say that
the solution
$y$
in
$\tilde{M}$
\textsl{is a normal pas-lifting} of the solution
$x$
\textsl{in}
$M,$
or, equivalently, that
$y$
is a \textsl{lifted pas-solution} of
$x,$
and also that
$x$
is the \textsl{pas-projection}
of
$y.$\\

\noindent We can use the previous remarks in combination with theorem
\ref{liftheorem} and we obtain the
following theorem
\begin{theorem}\label{liftheorem1}
Let $x(t),$ $t \in [\alpha, \beta]$ be an as-solution to
$(a, f).$ Then there is a lifted normal pas-solution
$y(t),$
$t \in [\alpha, \beta]$
to
$(\tilde{a}, \tilde{f})$
in
$\tilde{M}$
of a refinement of
$x(t).$
\end{theorem}
\paragraph{Equivalence between an IDE
$(a, f)$ in $M$ and
$(\tilde{a}_2, \tilde{f}_2)$ in $\tilde{M}_2.$}
Using the notation introduced in section \ref{sec-usingdesing}
we have obviously
that
$\tilde{M}_2$
is an open subset
of
$\tilde{M},$
$\tilde{\pi}_2 =\tilde{\pi}\vert \tilde{M}_2,$
$(\tilde{a}_2, \tilde{f}_2) = (\tilde{a}, \tilde{f})\vert \tilde{M}_2.$\\

\noindent
Let
$y(t),$
$t \in [\alpha, \beta],$
be a normal as-solution
to
$(a_k, f_k)$
in
$M^k$
and let
$\alpha \leq t_1
<....<t_r \leq \beta,$
be the set of all
$t \in [\alpha, \beta]$
such that
$x(t_i) \in M^k_0,$
as explained before the theorem.
Then we obtain a well defined collection of as-solutions
to
$(a_k, f_k)\vert M^k_2$
in
$M^k_2,$
namely,
$y\vert (t_i, t_{i + 1}),$
$i = 1,...,r-1,$
$y\vert (\alpha, t_1),$
if
$\alpha < t_1,$
$y\vert (t_r, \beta),$
if
$t_r < \beta.$
We will call this collection of solutions the \textsl{collection of solutions in
$M^k_2,$
or, also,
in
$\tilde{M}_2,$
induced by the as-solution}
$y.$
It is clear that by continuous extension, using lemma
\ref{lemmasuban1},
of each
$y\vert (t_i, t_{i + 1}),$
$i = 1,...,r-1,$
$y\vert (\alpha, t_1),$
if
$\alpha < t_1,$
$y\vert (t_r, \beta),$
if
$t_r < \beta,$
we obtain a refinement of
$x,$
with as-piece closures
$y\vert [t_i^+, t_{i + 1}^-],$
$i = 1,...,r-1,$
$y\vert [\alpha^+, t_1^-],$
if
$\alpha < t_1,$
$y\vert [t_r^+, \beta^-],$
if
$t_r < \beta.$
Of course, by continuity of
$y$
we have
$y\vert [t_i^+, t_{i + 1}^-] = y\vert [t_i, t_{i + 1}],$
$i = 1,...,r-1,$
$y\vert [\alpha^+, t_1^-] = y\vert [\alpha, t_1],$
if
$\alpha < t_1,$
$y\vert [t_r^+, \beta^-] = y\vert [t_r, \beta],$
if
$t_r < \beta.$\\

\noindent
Given any normal pas-solution
$y$
to
$(a_k, f_k)$
in
$M^k,$
we shall define the \textsl{collection of as-solutions in
$M^k_2,$
or, also, in
$\tilde{M}_2,$
induced by the normal pas-solution}
$y$
as being the union of the collections of as-solutions induced by all the as-pieces of
$y.$
Finally, given any normal pas-solution
$y$
to
$(\tilde{a}, \tilde{f}),$
in
$\tilde{M}.$
we shall define the \textsl{collection of as-solutions in
$\tilde{M}_2$
induced by the normal pas-solution}
$y$
as being the union of the collections of as-solutions
induced by all the as-pieces of which
$y$
is ultimately composed.
\\

\noindent
It is clear that a given normal
pas-solution
$y(t),$
$t \in [\alpha, \beta]$
to
$(\tilde{a}, \tilde{f})$
in
$\tilde{M},$
say
$y(t) \in M^k,$
for some
$k = 0,...,q$
as before,
cannot be completely recovered from its induced collection of as-solutions to
$(\tilde{a}_2, \tilde{f}_2)$
in
$\tilde{M}_2.$
This is because even if it is known that the value of
$y(t)$
at some
$t \in [\alpha, \beta]$
is the left or the right limit
at
$t$
of some of the as-pieces of the induced collection of as-curves by
$y,$
it is not known which one of those limits, unless one has some additional information,
like continuity
of
$y$
at
$t.$
Modulo this loss of information, we can consider
the IDE
$(\tilde{a}_2, \tilde{f}_2)$
with domain
$\tilde{M}_2$
and range
$F$
and
$(\tilde{a}, \tilde{f})$
with domain
$\tilde{M}$
and range
$F$
as being equivalent.
In this sense, the meaning of theorems \ref{liftheorem} and \ref{liftheorem1} is that solving a given analytic IDE $(a,
f)$
with domain
$M$
and
range
$F$
can be reduced essentially to solving a finite collection of
IDE of constant rank, namely, $(a_i, f_i)\vert M^i_2,$
$i = 0,...,q-1,$
or, equivalently, the system
of constant rank
$(\tilde{a}_2, \tilde{f}_2)$
with domain
$\tilde{M}_2$
and
range
$F.$
As we have said before it is clear that each $(a_i, f_i)\vert M^i_2,$ $i =
0,...,q-1$ defines an analytic family of vector fields and hence an analytic
affine distribution of constant rank
on the manifold $M^i_2.$
This immediately implies that
$(\tilde{a}_2, \tilde{f}_2)$ defines an analytic family of vector fields and hence an analytic
affine distribution of constant rank
on the manifold $\tilde{M}_2.$
Therefore, it also gives rise to an analytic control system,
~\cite{suss98:control}.
The usefulness of this kind of result relies obviously
on the interest that such control systems have in several fields.
For instance, the well known theorem of Sussmann,
~\cite{sussmann73}, tells us how to deal with problems of
reachability and observability, once one has a control system. See
also ~\cite{suss98:control} where this kind of questions are
solved in the context of subanalytic sets.
\paragraph{Reparametrization and extension of solutions.}
First we shall define the notion of a \textsl{reparametrization}
in the context of as-curves.
A \textsl{reparametrization}
of an as-curve
$x(t), t \in [t_0, t_1),$
in
$M,$
is a change of variables
$t = \tau(s),\, s \in [s_0, s_1),$
$(s \in (s_1, s_0]$)
such that
$\tau : [s_0, s_1) \rightarrow \Bbb{R}$
(respectively, $\tau : (s_1, s_0] \rightarrow \Bbb{R}$)
is an as-curve in
$\Bbb{R}$
which is also an homeomorphism
onto
$[t_0, t_1)$
such that
$\tau (s_0) = t_0$
and
$\tau \vert (s_0, s_1) : (s_0, s_1) \rightarrow (t_0, t_1)$
(respectively, $\tau : (s_1, s_0) \rightarrow (t_0, t_1)$)
is an analytic diffomorphism. It is easy to prove that in this case the composition
$(x \circ \tau)(s),$
$s \in [s_0, s_1)$
(respectively, $s \in (s_1, s_0]$)
is an as-curve in $M.$
Similar definitions and results hold for the case of as-curves
$x(t)$
whose domain is an interval of the type
$(t_0, t_1],$
$(t_0, t_1),$
$[t_0, t_1].$
\\

\noindent
Let
$(a, f)$
be a given IDE with domain
$M$
and range
$F.$
By definition,
$(a, f)$
is \textsl{homogeneous} if
$f = 0.$
It is clear that if
$x(t),\,\, t \in [t_0, t_1)$
is a given as-solution to an homogeneous system
$(a, 0)$
and
$t = \tau(s),\, s\in [s_0, s_1)$
($s \in (s_1, s_0]$)
is a given reparametrization
then the curve
$y(s) \equiv x \circ \tau (s), \, s\in [s_0, s_1)$
(respectively, $s \in (s_1, s_0]$)
is also an as-solution to
$(a, 0).$
More generally, if
$(a, f)$
is not necessarily homogeneous then $y(s)$ satisfies
$
a\left( y(s)\right)\dot{y}(s) = (d\tau/ds)f\left(y(s)\right), \, s\in (s_0, s_1),
$
or, using a different and also standard notation,
$
a\left( y(s)\right)\dot{y}(s) = \dot{t}(s)f\left(y(s)\right), \, s\in (s_0, s_1).
$
Similar results hold for the case of as-curves whose domain is an interval of the type
$(t_0, t_1],$
$(t_0, t_1),$
or
$[t_0, t_1].$
\begin{theorem}\label{thmas3ff}
$(a)$
Let
$(a, 0)$
be a given homogeneous IDE with domain
$M$
and range
$F.$
Let $x_{+}(t),$ $t\in [t_0, t_1),$ be an as-solution to
$(a, 0)$
in
$M,$
which is not a constant.
Then there is an as-solution $z(s),$ $s \in (s_0 - \delta_1, s_0 + \delta_2),$ in $M,$
for some $\delta_1, \delta_2 > 0,$
satisfying all the conditions stated in lemma \ref{thmas3}.
More precisely, $z(s_0) =
x_{+}(t_0),$ $z([s_0, s_0 + \delta_2)) = x_{+}([t_0, t_{2+}))$ for some $t_{2+} \in (0,
t_1)$ and we also have that $t_{2+}$ and $s_0 + \delta_2$ are such that
$x_{+}\vert
[t_0, t_{2+})$ and $z\vert [s_0, s_0 + \delta_2)$
are
homeomorphisms onto
$z([s_0, s_0 + \delta_2)) = x_{+}([t_0, t_{2+})),$
$x_{+}\vert
(t_0, t_{2+})$ and $z\vert (s_0, s_0 + \delta_2)$
are analytic diffeomorphisms onto
$z((s_0, s_0 + \delta_2)) = x_{+}((t_0, t_{2+}))$
which is an analytic submanifold which is also a semianalytic subset and, moreover,
the curve in
$\Bbb{R},$
say
$t = t(s),$
given by
$\left(x_{+}\vert (t_0, t_{2+})\right)^{-1} \circ (z\vert (s_0, s_0 + \delta_2)) : (s_0, s_0 + \delta_2) \rightarrow
\Bbb{R}$ is an as-curve in $\Bbb{R}$ which is an analytic
diffeomorfism onto its image $(t_0, t_{2+}).$ Moreover,
$t_{2+}$
can be chosen such that for each
$t_{3+} \in (t_0, t_{2+}],$
$x_{+}\left((t_0, t_{3+})\right)$ is a
neighborhood of $x_{+}(t_0)$ in $x_{+}\left([t_0, t_{2+})\right) - \{x_{+}(t_0)\}.$
Similar results hold for as-solutions
$x_{-}(t),$
$t\in (t_0, t_1].$
\\
\noindent
$(b)$
Let
$(a, f)$
be a given IDE with domain
$M$
and range
$F.$
Let $x_{+}(t),$ $t\in [t_0, t_1),$ be an as-solution to
$(a, f)$
in
$M$
which is not a constant.
Then there is an as-curve $z(s),$ $s \in (s_0 - \delta_1, s_0 + \delta_2),$ in $M,$
for some $\delta_1, \delta_2 > 0,$
not necessarily a solution, satisfying all the conditions stated in lemma \ref{thmas3},
as we have explained in $(a),$ and, besides,
the as-curve
$z(s), \,s\in [s_0,, s_0 + \delta_2)$
satisfies the equation
$a\left(z(s)\right)\dot z(s) = \dot{t}(s) f\left(z(s)\right),\, s\in (s_0, s_0 + \delta_2),$
with
$\dot{t}(s) > 0,$
$s \in (s_0, s_1).$
Similar results hold for as-solutions of the type
$x_{-} (t),$
$t\in (t_0, t_1].$
\\
\noindent
$(c)$
Let
$(a, f)$
be a given IDE with domain
$M$
and range
$F.$
Let $x_{+}(t),$ $t\in [t_0, t_1),$ be an as-solution to
$(a, f)$
in
$M$
which is not a constant
and let
$z(s)$
be as in
$(b).$
Then by conveniently diminishing the value of
$\delta_1$
and the value of
$\delta_2$
if necessary, we have the following. There is an as-solution
$x_{-}(t), \,t\in (t_{2-}, t_0]$
to some of the systems
$(a, \pm f)$
in
$M$
such that
$x_{-}(t_0) = x_{+}(t_0) = z(s_0),$
and
$x_{-}\vert
(t_{2-}, t_0]$
and
$z\vert (s_0-\delta_1, s_0]$ are homeomorphisms onto
$z((s_0 - \delta_1, s_0]) = x_{-}((t_{2-}, t_0]),$
$x_{-}\vert
(t_{2-}, t_0)$
and
$z\vert (s_0-\delta_1, s_0)$ are analytic diffeomorphisms onto
$z((s_0 - \delta_1, s_0)) = x_{-}((t_{2-}, t_0)),$
which is an analytic submanifold which is also a semianalytic subset.
Moreover,
$(x_{-}\vert (t_{2-}, t_0))^{-1} \circ (z\vert
(s_0 - \delta_1, s_0)): (s_0 - \delta_1, s_0) \rightarrow
\Bbb{R}$ is an as-curve in $\Bbb{R}$ which is an
as-diffeomorfism onto its image $(t_{2-}, t_0).$
We also have that
$t_{2-}$
can be chosen such that for each
$t_{3-} \in (t_{2-}, t_0],$
$x_{-}\left((t_{3-}, t_0)\right)$ is a
neighborhood of $x_{-}(t_0)$ in $x_{-}\left((t_{2-}, t_0]\right) - \{x_{-}(t_0)\}.$
The as-curve
$z(s), \,s\in [s_0, s_0 + \delta_2)$
satisfies the equation
$a\left(z(s)\right)\dot z(s) = \dot{t}(s)f\left(z(s)\right),\, s\in (s_0, s_0 + \delta_2)$
where
$\dot{t}(s) > 0,$ $s\in (s_0, s_0 + \delta_2).$
The as-curve
$z(s), \,s\in (s_0 - \delta_1, s_0]$
satisfies the equation
$a\left(z(s)\right)\dot z(s) = \dot{t}(s)f\left(z(s)\right),\, s\in (s_0 - \delta_1, s_0)$
where
$\dot{t}(s) > 0$ if
$x_{-}(t)$ satisfies
$(a, f)$
and
$\dot{t}(s) < 0$ if
$x_{-}(t)$ satisfies
$(a, -f).$
Similar results hold for as-solutions
$x_{-}(t),$ $t\in (t_0, t_1].$
\end{theorem}
\textbf{Proof.}
Part $(a)$ and part $(b)$ are a direct consequence
of lemma \ref{thmas3}.
To prove part $(c)$ consider the
homogeneous IDE $(\alpha, 0)$ with domain $\Bbb{R} \times M$ and
range $F$ where $\alpha (t, x)(\dot{t}, \dot{x}) = a(x)\dot{x} -
\dot{t}f(x).$ Consider the solution $(t(s), z(s)),\, s \in (s_0,
s_0 + \delta_2)$ to the homogeneous system $(\alpha, 0)$ where
$z(s),\, s \in (s_0 - \delta_1, s_0 + \delta_2)$ is the as-curve
considered in $(b).$
Since we know from $(b)$ that
$a\left(z(s)\right)\dot{z}(s)$ and $f\left(z(s)\right))$ are
linearly dependent for $s \in (s_0, s_0 + \delta_2)$ we can conclude, using
the analyticity of $z(s),\, s \in (s_0 - \delta_1, s_0 + \delta_2),$
that they must also be linearly dependent
for $s \in (s_0 - \delta_1, s_0 + \delta_2).$
We can also show, by using lemma \ref{lemmaas2},
that there exist
$\delta_1$
such that
$z\left((s_0 - \delta_1, s_0]\right) - \{z(s_0)\}$ is
nonempty and locally connected at $z(s_0),$
more precisely,
$z\left((s_0 - \delta_3, s_0)\right)$
is a neighborhood of
$z(s_0)$
in
$z\left((s_0 - \delta_1, s_0]\right) - \{z(s_0)\}$
for all
$s_0 - \delta_3 \in (s_0 - \delta_1, s_0).$
Moreover,
$\delta_1$
can be chosen such that
$z \vert [s_0 - \delta_1, s_0] : [s_0 - \delta_1, s_0] \rightarrow z \left([s_0 - \delta_1, s_0]\right)$
is an homeomorphism, $z \left((s_0 - \delta_1, s_0)\right)$
is an analytic submanifold which is a semianalytic
subset of $M$ and $z : (s_0 - \delta_1, s_0) \rightarrow z
\left((s_0 - \delta_1, s_0)\right)$ is an analytic
diffeomorphism.
Since we have a linear dependence between
$a\left(z(s)\right)\dot z(s)$ and $f\left(z(s)\right),\,
s\in (s_0 - \delta_1, s_0 + \delta_2)$
we have several cases.
Assume first that
$a\left(z(s)\right)\dot{z}(s) = 0,\, s\in (s_0, s_0 + \delta_2).$
Then since
$\dot{t}(s) \neq 0,\,s\in (s_0, s_0 + \delta_2)$
we can conclude that
$f\left(z(s)\right) = 0,\, s\in (s_0, s_0 + \delta_2).$
By analyticity of
$z(s),\, s\in (s_0 - \delta_1, s_0 + \delta_2)$
we obtain that
$a\left(z(s)\right)\dot{z}(s) = 0$ and $f\left(z(s)\right) = 0, \,s\in (s_0 - \delta_1, s_0 + \delta_2).$
Then to prove
$(c)$
in this case we can simply take
$t = s - s_0 + t_0,$
$t_{2-} = -\delta_1 + t_0$
and
$x_{2-}(t) = z(t + s_0 - t_0).$
We can proceed in a similar way if
we assume that
$f\left(z(s)\right) = 0,\,s\in (s_0, s_0 + \delta_2).$
Let us consider now the case where
$a\left(z(s)\right)\dot{z}(s),$
$s \in (s_0, s_0 + \delta_2)$
is not identically $0.$
We are going to show that,
after conveniently diminishing the value of
$\delta_1$ if necessary, there exists a unique as-curve in
$\Bbb{R},$ $\lambda :
(s_0 - \delta_1, s_0) \rightarrow \Bbb{R}$ such that
$a\left(z(s)\right)\dot z(s) = \lambda(s)f\left(z(s)\right),\,
s\in (s_0 - \delta_1, s_0).$
First one should take into
account that $z(s),\, s\in (s_0 - \delta_1, s_0 + \delta_2)$ is
analytic then so are
$f\left(z(s)\right)$
and
$a\left(z(s)\right)\dot{z}(s),\, s\in (s_0 - \delta_1, s_0 +
\delta_2)$ and therefore they have at most a finite number of
isolated zeros in a neighborhood of
$s_0.$
It is easy to see from the equation
$a\left(z(s)\right)\dot z(s) = \lambda(s)f\left(z(s)\right),\,
s\in (s_0 - \delta_1, s_0 + \delta_2)$
that
$\lambda(s),$
or rather its extension for complex
$s,$
is a meromorphic function in a neighborhood of
$s_0.$
On the other hand, since
$\lambda(s) = \dot{t}(s),$
$s \in (s_0, s_0 + \delta_1)$
where
$t(s)$
is bounded, we have that
$\lambda(s)$
cannot have a pole
at
$s_0$
therefore it must be analytic in a neighborhood of
$s_0.$
This implies that by conveniently diminishing the value of
$\delta_1$
and the value of
$\delta_2$
if necessary,
we have that
$\lambda(s),$
$s \in (s_0 - \delta_1, s_0 + \delta_2)$
is real analytic.
By
conveniently diminishing the value of $\delta_1$ if necessary we
can assume without loss of generality that $\lambda(s) \neq 0, \,
s\in (s_0 - \delta_1, s_0).$
Let us assume first that
$\lambda > 0.$
Then the result follows by taking $t
= t(s),\, s\in (s_0 - \delta_1, s_0 + \delta_2)$ such that $dt/ds =
\lambda(s),$ $t(s_0) = t_0$ and $t(s_0 - \delta_1) = t_{2-},$
which defines
$t(s),\, s\in (s_0 - \delta_1, s_0]$
as an as-curve in
$\Bbb{R}$
and also
$t_{2-}.$
Since
$\lambda(s) > 0, \,
s\in (s_0 - \delta_1, s_0)$
we have that
$t(s)$
is an analytic diffeomorphism
from
$(s_0 - \delta_1, s_0)$
onto
$(t_{2-}, t_0)$
which is an as-curve,
while
$t : [s_0 - \delta_1, s_0] \rightarrow [t_{2-}, t_0]$
is an homeomorphism.
Then we can define
$x_{2-}(t),$
$t \in (t_{2-}, t_0]$
by
$x_{2-}(t) = z\left(s(t)\right)$
where
$s(t)$
is the inverse of
$t(s).$
The rest of the proof can be performed in a similar way.
\quad $\blacksquare$\\

\noindent
The previous theorem says, in particular,
that an as-solution $x_{+}(t),\, t\in [t_0, t_1)$ to a given IDE
$(a, f)$ gives rise to an as-solution $x_{-}(t),\, t\in
(t_{2-}, t_0]$
of
$(a, f)$
or
$(a, -f)$
such that $x_{+}(t_0) = x_{-}(t_0)$
and
$\operatorname{graph}x_{-} \cup \operatorname{graph}x_{+}
=
\operatorname{graph}z$
for some as-curve
$z(s),\, s\in (s_0 - \delta_1, s_0 + \delta_2),$
which is a solution
to the homogeneous system
$a(z)\dot{z} = \dot{t}f(z),$
$s \in (s_0 - \delta_1, s_0) \cup (s_0, s_0 + \delta_2),$
and
$x_{-},$
$x_{+}$
are reparametrizations of
$z\vert (s_0 - \delta_1, s_0],$
$z\vert [s_0, s_0 + \delta_2).$\\

\noindent
For a given IDE
$(a, f)$
with domain
$M$
and
range
$F,$
an \textsl{as-solution}
$x(t),$ $t \in [t_0, t_1)$ ($t \in (t_0, t_1]$, $t \in (t_0, t_1)$, $t \in [t_0, t_1]$),
\textsl{to}
$(a, \pm f),$
in
$M$
is, by definition, an as-solution
to some of the systems
$(a, f)$
or
$(a, -f),$
in
$M.$
The notion of a pas-solution or
an lcs-solution
to
$(a, \pm f),$
in
$M$
is defined by the condition
that each as-piece is an as-solution to
$(a, \pm f)$
in
$M.$\\

\noindent
We have the following extension theorem
\begin{theorem}\label{extensionthm2}
Let
$(a, f)$
be an IDE with domain
$M$
and range
$F.$
Let
$x(t),\, t\in [t_0, t_1]$
be a continuous pas-solution
to
$(a, f),$
so its graph is compact.
Then $x(t)$ can be extended to a continuous pas-solution
$\bar{x}(t),\, t\in [\bar{t}_0, \bar{t}_1],$
to
$(a, \pm f)$
for some
$\bar{t}_0 < t_0 < t_1 < \bar{t}_1.$
\end{theorem}
\textbf{Proof.}
If
$x$
is a constant the proof is immediate. Let us assume that
$x$
is not a constant.
We have an as-piece
of
$x$
of the type
$x\vert [t_0, t_0 + \delta_2).$
As a consequence of theorem
\ref{thmas3ff}
we have that
$x\vert [t_0, t_0 + \delta_2)$
can be extended to a continuous pas-solution, say
$u(t),$
$t \in [t_{2-}, t_0 + \delta_2),$
to
$(a, \pm f),$
where
$u_{+} = x\vert [t_0, t_0 + \delta_2)$
and
$u_{-} : ([t_{2-}, t_0] \rightarrow M$
are the as-pieces.
We can glue the as-piece
$u_{-}$
to
$x$
and take
$\bar{t}_0 = t_{2-}$
to obtain a continuous
pas-solution
to
$(a, \pm f)$
defined
in
$[\bar{t}_0, t_1)$
which extends
$x.$
We can proceed now in a similar way to obtain, in turn, an extension of this solution to a pas-continuous solution to
$(a, \pm f)$
defined
in some interval
$[\bar{t}_0, \bar{t}_1].$
\quad $\blacksquare$\\
\begin{corollary}\label{extensionthm2}
Let
$(a, 0)$
be an homogeneous IDE with domain
$M$
and range
$F.$
Let
$x(t),\, t\in [t_0, t_1]$
be a continuous pas-solution
to
$(a, f),$
so its graph is compact.
Then $x(t)$ can be extended to a continuous pas-solution
$\bar{x}(t),\, t\in [\bar{t}_0, \bar{t}_1],$
to
$(a, 0)$
for some
$\bar{t}_0 < t_0 < t_1 < \bar{t}_1.$
\end{corollary}
\section{The example of the symmetric elastic sphere}\label{section-ERS}
\noindent As we have said before the main purpose of this paper is
not to pursue the investigation of how to obtain systematically
desingularizations of given IDE by using systematic procedures to
desingularize given closed analytic sets and then apply all this to solve a given  IDE
by a systematic procedure,
but rather to show how having a
desingularization of a given IDE may help to understand its
solutions. In the next section we study an example from
nonholonomic mechanics, namely the symmetric elastic sphere, and
show using a desingularization of the equations of motion reduced
by the symmetry, which constitute an IDE, how to study the
dynamics and solve completely the system by quadratures.
\paragraph{Preliminaries.}

Continuing with the description of the symmetric elastic sphere
given in section \ref{INTR} we shall now give a more precise description of this system as a nonholonomic system. For references on nonholonomic mechanics related to rolling spheres see
~\cite{MS, BKMM, CEMARA2, fufaev, CMDLI, CELARE} and references therein.\\

\noindent
A rigid sphere rolling on a plane can be
modelled as a nonholonomic system on the group $SO(3) \times \Bbb
R^2$ where, for a given element $(A, x) \in SO(3) \times \Bbb
R^2,$ $A$ represents a rigid rotation and $x$ the position of the point of contact of the sphere with the plane.
Now we will assume that
$SO(3) \times \Bbb R^2$
is also a good configuration space for an elastic sphere where
deformations are small and concentrated only near the area of contact which is a
small circle whose center has a position given by
$x.$
The kinematics of this system can be
conveniently described as follows. We assume that there is an
ortonormal system fixed in the space, say $(e_1, e_2, e_3),$ $e_1
= (1, 0, 0),$ $e_2 = (0, 1, 0),$ $e_3 = (0, 0, 1),$ then we have a
basis moving with the body, $(Ae_1, Ae_2, Ae_3),$ where $A =
A(t).$ We introduce the variable $z \in S^2,$ given by $z = A
e_3.$ The spatial angular velocity $\omega$ can be written $\omega
= v_0 z + z \times \dot{z},$ so $v_0 = \langle \omega, z \rangle$
is the component of $\omega$ along $z.$ The nonholonomic
constraint is given by the nonsliding condition. Since the area of
contact of the sphere with the plane is a circle we must add to
the usual nonsliding condition $\omega \times r e_3 = \dot{x}$ for
the rigid rolling sphere (see, for instance, ~\cite{CELARE} where Lagrange-D'Alembert-Poincar\'{e} equations for the rigid sphere have been written)
where
$r$
is the radius of the sphere,
the extra condition that the vertical
component of the spatial angular velocity is $0,$ that is,
$\omega_3 = 0.$\\

\noindent We are going to assume that the center of mass
coincides with the center of the sphere and that the principal
axis of inertia are $(Ae_1, Ae_2, Ae_3).$ The three principal moments of
inertia of the sphere are $I_1,$ $I_2,$ $I_3,$ and we are going to
assume that $I_1 = I_2.$ We introduce the adimensional
quantities $\alpha = I_3 / I_1$ and $\beta = Mr^2/ I_1,$ where $M$
is the mass of the sphere.
The Lagrangian of the system is given by the kinetic energy,
$$
\frac{1}{2}I_1 \dot{z}^2
+ \frac{1}{2}I_3 v_0^2
+ \frac{1}{2}M\dot{x}^2
$$
where
$\dot{x}$
is the velocity of the center of the sphere. The nonholonomic constraint is given by
$\dot{x} = \omega \times r e_3$
and
$\omega_3 = 0,$
and using this we can conclude that the kinetic energy of the actual motion of the symmetric sphere is given by

$$
E = \frac{1}{2}(I_1 + Mr^2)\dot{z}^2 + \frac{1}{2}(I_3 + Mr^2)v_0^2.
$$
\\

\noindent
As we have said in the Introduction the addition of the extra condition
$\omega_3 = 0$
introduces an extra singularity in the reduced system,
which
is an IDE, and we need to apply our desingularization procedure to
obtain a single equivalent differential equation describing the
system, in the sense of theorems \ref{liftheorem} and \ref{liftheorem1}. We will show
that the desingularizing manifold containing the essential
dynamics, in which this differential equation is defined, is
diffeomorphic to $S^2 \times S^1.$ Another interesting feature of
the desingularization method for the example of the symmetric
elastic sphere is that integrability by quadratures
appears in a natural way, in terms of angular coordinates of
$S^2 \times S^1.$\\

\paragraph{The IDE for the symmetric elastic sphere.}
\noindent As a result of reduction by the symmetry techniques, in
this case reduction by the subgroup $SO(2) \times \Bbb R^2$ we obtain the following system of
Lagrange-D'Alembert-Poincar\'{e} equations, which is an IDE,

\begin{eqnarray}
\label{esf1} \nonumber (\alpha + \beta)(z \times e_3)\dot{v}_0
+
(1+\beta)<z,e_3>\nabla_{\dot{z}}{\dot{z}}- & & \\
(\alpha+\beta)v_0 <z,e_3> (z\times\dot{z}) & = & 0 \\
\label{esf2} v_0<z,e_3>+<z \times \dot{z},e_3> & = & 0.
\end{eqnarray}
Here $\nabla$ represents the Levi-Civita connection on $S^2$ with
respect to the standard metric. This is a consequence of the
methods developed in ~\cite{CEMARA2}, after some more or less
straightforward calculations which we will not explain here.
The previous Lagrange-D'Alembert-Poincar\'{e} equations are derived
under the assumption $z_3 \neq 0$ because the so called
\textsl{dimension assumption} adopted in ~\cite{CEMARA2} is not
satisfied for the whole manifold $S^2.$
Nevertheless, by
continuity, the equations (\ref{esf1}), (\ref{esf2})
are also satisfied by the motion of the rolling ball for
$z_3 = 0.$ Without using the derivation of the
Lagrange-D'Alembert-Poincar\'{e} equations, the careful reader may
want to check directly that the previous system of equations, or equivalently, the system of equations
(\ref{ec1})-(\ref{ec77}),
is
equivalent to balance of momentum plus the condition
$\omega_3 = 0.$\\

\noindent Since $\nabla_{\dot{z}}{\dot{z}}=z \times (\ddot{z}
\times z)$ by taking the inner product of (\ref{esf1}) with
$\dot{z}$ and using (\ref{esf2}) we get, at least for
$<z, e_3> \neq 0,$
\begin{equation}
\label{energ}
0
=\frac{d}{dt}\left((1+\beta)\dot{z}^2+(\alpha+\beta)v_0^2\right),
\end{equation}
from which one deduces
\begin{equation}
\label{energ} 2\epsilon=(1+\beta)\dot{z}^2+(\alpha+\beta)v_0^2,
\end{equation}
where $\epsilon$ represents the normalized energy. This equation
represents conservation of energy, as one can check more directly
by looking at the expression of the kinetic energy
$E$
given at the beginning of this section. We shall assume from now on that $\epsilon > 0,$
otherwise the motion is trivial.\\

\noindent We have the following equations to be satisfied for
the symmetric elastic sphere in variables $(z,u)$ where
$\dot{z}=v$ and $v \times z=u,$ so the variable
$v_0$
does not appears,
\noindent
\begin{eqnarray}
\label{Esf2} (1+\beta)<z,e_3><\dot{u},e_3 \times z>+(\alpha +
\beta)<u,e_3>^2 & = & 0\\
\label{Esf1}
(1+\beta)<z,e_3>^2u^2+(\alpha+\beta)<u,e_3>^2-2\epsilon<z,e_3>^2 & = & 0,
\end{eqnarray}
Equation (\ref{Esf2}) is obtained by taking the inner product of (\ref{esf1}) with $e_3$ and
using (\ref{esf2}) while equation
(\ref{Esf1}) is obtained from equation (\ref{energ}) and equation (\ref{esf2}).
\\

\noindent Equations (\ref{Esf2}) and (\ref{Esf1}) involve the variables
$(z, u) \in TS^2$ and we have a natural inclusion
$TS^2 \subseteq S^2 \times \Bbb R^3$, where
\[
TS^2=\{(z,u) \in S^2 \times \Bbb R^3\,:\,<z,u>=0\}.
\]

\noindent Equations (\ref{Esf2}), (\ref{Esf1}) form
an IDE in the manifold $TS^2.$ By adding the equation $z ^2 = 1,$
we obtain an equivalent IDE in the variables $(z, u) \in \Bbb R^3
\times \Bbb R^3.$ Let us include, in addition, the equation of
conservation of energy (\ref{energ}), written in terms of the
variables $(z, u, v_0),$ as follows
\begin{equation}
\label{energ22} 2\epsilon=(1+\beta)u^2+(\alpha+\beta)v_0^2.
\end{equation}
In other words, we are going to study the system of equations
given by (\ref{Esf2}), (\ref{Esf1}), (\ref{energ22})
for a fixed $\epsilon > 0.$
Of course equations (\ref{Esf1}) and
(\ref{energ22}) taking into account (\ref{esf2})
are redundant for $z_3 \neq 0,$ but for $z_3 = 0$
the system given by the equations
(\ref{Esf2}), (\ref{Esf1}) includes solutions of the type $z(t) = const,$ with
$z_3 = 0,$ where $v_0$ takes any given value and, since $u = 0$
for this kind of motion, the energy is given by $2\epsilon =
(\alpha + \beta)(v_0)^2$ and therefore the condition that the
energy must have a fixed value will not be satisfied.
Of course we can study with our methods both, the system given by
(\ref{Esf2}), (\ref{Esf1}), (\ref{energ22}) and also
the system given by (\ref{Esf2}), (\ref{Esf1}), but
we will chose to study just the first of them, for simplicity.

\paragraph{The IDE for the symmetric elastic sphere in the standard form.}
Considering that $z^2=1$ and $\dot{z}=z \times u,$ we must have
$<z, u> = 0.$ Using what was said in the previous paragraph, and
according to theorem \ref{thm0} we can write the system of
equations for the symmetric elastic sphere in variables
$(z,u,v_0) \in \Bbb R^3 \times \Bbb R^3\times \Bbb R$ as follows,
\begin{eqnarray}
\label{ec1}
\dot{z}_1 & = & z_2u_3-z_3u_2 \\
\label{ec2}
\dot{z}_2 & = & z_3u_1-z_1u_3 \\
\label{ec3}
\dot{z}_3 & = & z_1u_2-z_2u_1\\
\label{ec4}
0 & = & (1+\beta)z_3(-z_2\dot{u}_1+z_1\dot{u}_2)+(\alpha + \beta)u_3^2\\
\label{ec5} 0 & = & (1+\beta)z_3^2(u_1^2+u_2^2+u_3^2)+(\alpha
+ \beta)u_3^2-2\epsilon z_3^2\\
\label{ec6}
0 & = & z_1^2+z_2^2+z_3^2 -1\\
\label{ec7}
0 & = & z_1 u_1 + z_2 u_2 + z_3 u_3\\
\label{ec77} 0 & = & 2\epsilon - (1+\beta)u^2-(\alpha+\beta)v_0^2.
\end{eqnarray}
The system (\ref{ec1})-(\ref{ec77}) can be written in the form
\[a(X)\dot{X} = f(X),\]
with $X=(z,u, v_0),$ where,
\[a(z,u, v_0)
= \left [\begin{array}{ccccccc}
1 & 0 & 0 & 0 & 0 & 0 & 0\\
0 & 1 & 0 & 0 & 0 & 0 & 0 \\
0 & 0 & 1 & 0 & 0 & 0 & 0 \\
0 & 0 & 0 & -(1+\beta)z_2z_3 & (1+\beta)z_1z_3 & 0 & 0 \\
0 & 0 & 0 & 0 & 0 & 0 & 0\\
0 & 0 & 0 & 0 & 0 & 0 & 0\\
0 & 0 & 0 & 0 & 0 & 0 & 0\\
0 & 0 & 0 & 0 & 0 & 0 & 0\\
\end{array}
\right];\]
\[
 f(z,u, v_0)
 =
 \left
 [\begin{array}{c}z_2u_3-z_3u_2 \\
              z_3u_1-z_1u_3 \\
              z_1u_2-z_2u_1\\
              -(\alpha+\beta)u_3^2\\
              (1+\beta)z_3^2(u_1^2+u_2^2+u_3^2)+(\alpha
              + \beta)u_3^2-2\epsilon z_3^2\\
z_1^2+z_2^2+z_3^2-1\\
z_1 u_1 + z_2 u_2 + z_3 u_3\\
2\epsilon - (1+\beta)u^2-(\alpha+\beta)v_0^2
\end{array}\right].\]

\paragraph{Application of the algorithm.}
We will work on the manifold $M = \Bbb R ^7,$ where $(z_1, z_2,
z_3, u_1, u_2, u_3, v_0) \in \Bbb R ^7$ are independent variables.
Then our IDE is given by equations (\ref{ec1})-(\ref{ec77}). We
can easily see that $k_r = 4,$ $S_4(M) = M,$ $L_4(M) = M_0,$ $M_1
= M-L_4(M),$ $M_2 = \emptyset.$ Now we shall describe $M_0$ by
equations. Let
\begin{align}
\varphi_1 &=
-(1 + \beta) z_2 z_3\\
\varphi_2 &=
(1 + \beta) z_1 z_3\\
\nu_1 &=
(1+\beta)z_3^2(u_1^2+u_2^2+u_3^2)+(\alpha + \beta)u_3^2-2\epsilon z_3^2\\
\nu_2 &=
z_1^2+z_2^2+z_3^2-1\\
\nu_3 &=
z_1 u_1 + z_2 u_2 + z_3 u_3\\
\nu_4 &= 2\epsilon - (1+\beta)u^2-(\alpha+\beta)v_0^2.
\end{align}
As we know $M_0 = L_4(M)$ is given by the condition that
$\operatorname{rank}[a, f] \leq 4.$ Let
\begin{align}
M_{0a} &=
\{\varphi_1 = 0, \varphi_2 = 0\}\\
&=
\{z_3 = 0\} \cup \{z_1 = 0, z_2 = 0\}\\
M_{0b} &= \{\nu_1 = 0, \nu_2 = 0, \nu_3 = 0, \nu_4 = 0\},.
\end{align}

\noindent Then we can easily see that $M_0 = M_{0a} \cup M_{0b}.$
The desingularization $M^1$ of $M_0$ will be the disjoint union of
the desingularizations of $M_{0a}$ and
$M_{0b}.$\\

\noindent
The desingularization $M_a^1$ of $M_{0a}$ can be
described by $M_a^1 \equiv \{z_3 = 0\} \bigsqcup \{z_1 = 0, z_2 =
0\},$ where $\bigsqcup$ means \textsl{disjoint union} and the
projection $\pi_0$ is the identity on each disjoint piece of $M^1_a.$ One can see using
(\ref{ec1}) - (\ref{ec77}) that the lifted system $(a_1, f_1)\vert
\{z_3 = 0\}$ satisfies $z_3 = 0,$ $u_3 = 0,$ $z_1 ^2 + z_2^2 = 1,$
which implies $\dot{z} = 0,$ and also, since $u = \dot{z} \times
z,$ that $u = 0.$ This describes the motion completely. It
consists of the rolling of the sphere with $z(t) = (z_{10},
z_{20}, 0)$ fixed and the $z$ component of the angular velocity
$v_0$ satisfies $2\epsilon = (\alpha + \beta)(v_0)^2.$ The lifted
system $(a_1, f_1)\vert \{z_1 = 0, z_2 = 0\}$ satisfies $z_1 = 0,$
$z_2 = 0,$ $z_3 = \pm 1,$ therefore $\dot{z} = 0,$ and then $u =
0,$ which contradicts equation $\nu_1 = 0,$ because we have
assumed $\epsilon > 0.$ So there is no motion, that is, no
solution, for the system $(a_1,
f_1)\vert \{z_1 = 0, z_2 = 0\}.$\\

\noindent
 Now we will desingularize $M_{0b}.$ We are going to see
that $M_{0b}$ is in fact a nonsingular manifold. More precisely, we will
define the desingularizing manifold $M^1_b$ by equations in the
variables $(z,u,v_0)$, with $v_0 z_3 = u_3$, from (\ref{esf2}).
For simplicity, we call $\mu = 2\epsilon / (1+\beta) > 0$ and
$\lambda = (\alpha + \beta )/(1+\beta) > 0,$ from now on. Then we
have the following equations defining the nonsingular manifold
$M^1_b,$

\begin{eqnarray}
\label{EEc5}
0 & = & u_3-v_0z_3\\
\label{EEc6}
0 & = & u_1^2+u_2^2+u_3^2+ \lambda v_0^2-\mu \\
\label{EEc7}
0 & = & z_1^2+z_2^2+z_3^2-1 \\
\label{EEc8} 0 & = & z_1u_1+z_2u_2+z_3u_3.
\end{eqnarray}
See {\bf Appendix A} for a proof that the system above defines a
nonsingular manifold. The map $\pi_0 : M^1_b \rightarrow M$ is
then given by the restriction of the identity $(z, u, v_0)
\rightarrow (z, u, v_0)$ to $M^1_b$ and one can check that the
image of $\pi_0$ is precisely
$M_{0b}.$\\

\noindent According to theorem \ref{thm0}, the system lifted to
$M^1_b$ has the same solutions as the system given by the
equations

\begin{eqnarray}
\label{Ec1}
\dot{z}_1 & = & z_2u_3-z_3u_2 \\
\label{Ec2}
\dot{z}_2 & = & z_3u_1-z_1u_3 \\
\label{Ec3}
\dot{z}_3 & = & z_1u_2-z_2u_1 \\
\label{Ec4}
z_2\dot{u}_1-z_1\dot{u}_2 & = & \lambda v_0u_3 \\
\label{Ec5}
0 & = & u_3-v_0z_3\\
\label{Ec6}
0 & = & u_1^2+u_2^2+u_3^2+ \lambda v_0^2-\mu \\
\label{Ec7}
0 & = & z_1^2+z_2^2+z_3^2-1 \\
\label{Ec8} 0 & = & z_1u_1+z_2u_2+z_3u_3.
\end{eqnarray}
\\

\noindent More precisely, if we define
\[
\tilde{a}(z,u,v_0)=\left[\begin{array}{ccccccc}
1 & 0 & 0 & 0 & 0 & 0 & 0\\
0 & 1 & 0 & 0 & 0 & 0 & 0\\
0 & 0 & 1 & 0 & 0 & 0 & 0\\
0 & 0 & 0 & z_2 & -z_1 & 0 & 0 \\
0 & 0 & 0 & 0 & 0 & 0 & 0\\
0 & 0 & 0 & 0 & 0 & 0 & 0\\
0 & 0 & 0 & 0 & 0 & 0 & 0\\
0 & 0 & 0 & 0 & 0 & 0 & 0
\end{array}\right], \tilde{f}(z,u,v_0)=\left[\begin{array}{c}

z_2u_3-z_3u_2 \\
z_3u_1-z_1u_3 \\
z_1u_2-z_2u_1 \\
\lambda v_0u_3 \\
u_3-v_0z_3\\
u_1^2+u_2^2+u_3^2+\lambda v_0^2-\mu \\
z_1^2+z_2^2+z_3^2-1\\
z_1u_1+z_2u_2+z_3u_3 \end{array} \right],
  \]
we see that the system (\ref{Ec1})-(\ref{Ec8}) is in the form
$\tilde{a}(y)\dot{y}= \tilde{f}(y),$ with $y=(z,\,u,\,v_0)$, so it
is an IDE in standard form with domain $\Bbb{R}^7$ and range
$\Bbb{R}^8$, and our IDE in $M^1_b$ is given by the restriction
$(a_1, f_1) = (\tilde{a}, \tilde{f})\vert M^1_b.$
In order to continue with the algorithm, we shall find explicitly
the lifted system $(a_1, f_1).$
By differentiating the equations (\ref{EEc5})-(\ref{EEc8}),
eliminating the redundant equation $z_1 \dot{z}_1 + z_2 \dot{z}_2
+ z_3 \dot{z}_3 = 0,$ also realizing appropriate linear operations
in the range space, that is, using projections as explained in
section \ref{sec-difalgsys}, and also according to theorem
\ref{thm0}, we have the following system with domain $\Bbb{R}^7$
and range $\Bbb{R}^{11},$ which is also equivalent to our system
$(a_1, f_1),$
\begin{eqnarray}
\label{Ec11}
\dot{z}_1 & = & z_2u_3-z_3u_2 \\
\label{Ec22}
\dot{z}_2 & = & z_3u_1-z_1u_3 \\
\label{Ec33}
\dot{z}_3 & = & z_1u_2-z_2u_1 \\
\label{Ec44}
z_2\dot{u}_1-z_1\dot{u}_2 & = & \lambda v_0u_3 \\
\label{Ec55}
z_1\dot{u}_1+z_2\dot{u}_2+z_3\dot{u}_3 & = & 0 \\
\label{Ec66}
u_1\dot{u}_1+u_2\dot{u}_2+u_3\dot{u}_3+\lambda v_0\dot{v}_0 & = & 0\\
\label{Ec77}
\dot{u}_3-z_3\dot{v}_0 & = & v_0z_1u_2-v_0z_2u_1\\
\label{Ec51}
0 & = & u_3-v_0z_3\\
\label{Ec61}
0 & = & u_1^2+u_2^2+u_3^2+ \lambda v_0^2-\mu \\
\label{Ec71}
0 & = & z_1^2+z_2^2+z_3^2-1 \\
\label{Ec81} 0 & = & z_1u_1+z_2u_2+z_3u_3.
\end{eqnarray}

\noindent This system is still not completely desingularized. One
can check by direct calculation that it can be desingularized in
two more iterations of the algorithm. However, in this example
there is an interesting alternative to find the solutions, which
starts with a precise description of the manifold $M^1_b.$ We
prefer this alternative because having an identification of
$M^1_b$ also helps to understand the dynamics in a direct way, as
we will see soon. This is precisely one of the points of the
present work, namely, to show how in examples the
desingularization method may helps to understand the dynamics of a
given system.
\paragraph{Identification of $M^1_b.$}
The manifold $M^1_b$ is given by the equations (\ref{EEc5}) -
(\ref{EEc8}) in the space of the variables
$(z_1,z_2,z_3,u_1,u_2,u_3,v_0)$,
as we have seen before.
The equation (\ref{EEc8}) tells us that $u$ is a vector
tangent to the $2$-sphere $S^2,$ given by $z^2-1=0.$
Heuristically, for each $z \in S^2$ we consider the 3 dimensional
space $T_{z}{S^2} \times R_z,$ where $R_z$ represents a line
normal to the sphere at $z \in S^2,$ so the collection of all
$R_z$ is a trivial real line vector bundle with base $S^2.$
Equation (\ref{EEc5}) is a plane containing the origin $0=0_z$
since $z_3$ is fixed once $z$ is fixed. Equation (\ref{EEc6})
gives an ellipsoid. The intersection of the plane with the
ellipsoid is an ellipse.
Therefore
$M^1_b$
must be some fiber bundle with fiber
$S^1$
and base
$S^2.$
Using all this and some imagination we can see that it is, in fact,
the trivial bundle
$S^2 \times S^1,$
moreover, we have the following parametrization of $M^1_b$ in variables
$(\theta,\,\varphi,\,\psi).$
In any case, this assertion can be easily checked after some
straightforward calculations.
\begin{eqnarray}
\label{par1}
z_1 & = & \sin{\theta}\cos{\varphi} \\
\label{par2}
z_2 & = & \sin{\theta}\sin{\varphi}\\
\label{par3}
z_3 & = & \cos{\theta}\\
\label{par4} u_1 & = &
-a\cos(\varphi-\psi)\cos^2{\theta}\cos{\varphi}
-b\sin(\varphi-\psi)\sin{\varphi}\\
\label{par5} u_2 & = &
-a\cos(\varphi-\psi)\cos^2{\theta}\sin{\varphi}
+b\sin(\varphi-\psi)\cos{\varphi}\\
\label{par6}
u_3 & = & a\cos(\varphi-\psi)\cos{\theta}\sin{\theta}\\
\label{par7} v_0 & = & a\cos(\varphi-\psi)\sin{\theta},
\end{eqnarray}
where
\[a=\sqrt{\frac{\mu}{\lambda\sin^2{\theta}
+\cos^2{\theta}}},\,\,\,b=\sqrt{\mu}'.\] In other words, by some
straightforward calculations we can check that
$(z_1,z_2,z_3,u_1,u_2,u_3,v_0)$ in coordinates
$(\theta,\,\varphi,\,\psi)$ satisfies (\ref{EEc5})-(\ref{EEc8}).\\

\noindent We can see that equations (\ref{par1})-(\ref{par7})
define a diffeomorphism $f:S^2 \times S^1 \rightarrow M^1_b$,
$f(z,(\cos{\psi},\sin{\psi}))=(z,u,v_0),$ which gives the desired
identification of $M^1_b.$ See {\bf Appendix B} for a proof that
$f$ is a diffeomorphism.
\\

\paragraph{The differential equation for
the symmetric elastic  sphere in $M^1_b$ in variables $(\theta,
\varphi, \psi).$} Considering the parametrization for $M^1_b$
given by (\ref{par1})-(\ref{par7}) we get the equations for
(\ref{Ec1})-(\ref{Ec8}) in coordinates $(\theta,\varphi,\psi)$,
\begin{eqnarray}
\label{pol1} \cos{\theta}\cos{\varphi}\dot{\theta}
-\sin{\theta}\sin{\varphi}\dot{\varphi}
& = & a\cos{\theta}\sin{\varphi}\cos(\varphi-\psi)-\\
\nonumber
&  & -b\cos{\theta}\cos{\varphi}\sin(\varphi-\psi)\\
\label{pol2} \cos{\theta}\sin{\varphi}\dot{\theta}
+\sin{\theta}\cos{\varphi}\dot{\varphi}
& = & -a\cos{\theta}\cos{\varphi}\cos(\varphi-\psi)-\\
\nonumber
& & -b\cos{\theta}\sin{\varphi}\sin(\varphi-\psi)\\
\label{pol3}
-\sin{\theta}\dot{\theta} & = & b\sin{\theta}\sin(\varphi-\psi) \\
\label{pol4}
a\sin{\theta}\cos^2{\theta}\cos(\varphi-\psi)\dot{\varphi}-\\
\nonumber
  -b\sin{\theta}\cos(\varphi-\psi)
(\dot{\varphi}-\dot{\psi}) & = & \lambda
a^2\cos^2(\varphi-\psi)\sin^2{\theta}\cos{\theta}
\end{eqnarray}

\noindent If $\sin{\theta}\not=0$ the system (\ref{pol1})
-(\ref{pol4}) becomes
\begin{eqnarray}
\label{Pol1}
\dot{\theta} & = & -b\sin(\varphi-\psi)\\
\label{Pol2} \dot{\varphi}
& = & -a\frac{\cos{\theta}}{\sin{\theta}}\cos(\varphi-\psi) \\
\label{Pol3} \dot{\psi} & = & a\cos(\varphi-\psi)
\frac{\cos{\theta}}{\sin{\theta}}\left(\frac{b}{a}-1\right),
\end{eqnarray}
or equivalently,
\begin{eqnarray}
\label{PPol1}
\dot{\theta} & = & -b\sin(\varphi-\psi)\\
\label{PPol2} \dot{\varphi}
& = & -a\frac{\cos{\theta}}{\sin{\theta}}\cos(\varphi-\psi) \\
\label{PPol3} \dot{\psi} & = &
(b-a)\frac{\cos{\theta}}{\sin{\theta}}\cos(\varphi-\psi)
\end{eqnarray}

\noindent It can be easily seen that this system can be integrated
by quadratures. For instance, if we call $w=\varphi-\psi,$ we can
write (\ref{PPol1})-(\ref{PPol3}) as a planar system in
coordinates $(\theta,w),$
\begin{eqnarray}
\label{planar1}
\dot{\theta} & = & -b\sin{w}\\
\label{planar2} \dot{w} & = &
-b\frac{\cos{\theta}}{\sin{\theta}}\cos{w},
\end{eqnarray}
which in turn leads to the separable equation
\begin{equation}
\label{escalar} \frac{d\theta}{dw}=\tan{\theta}\tan{w}.
\end{equation}

\noindent Of course the system (\ref{pol1}) - (\ref{pol4}) is
still an analytic IDE on the analytic
manifold $S^2 \times S^1$ which is of constant rank for $\sin
\theta \neq 0,$ and the rank changes for $\sin \theta = 0.$ So we
should continue the desingularization process,
which is not difficult but we are not going to explain the
details.
Instead, we simply observe that it is not difficult to see,
alternatively, that the only solution with some initial condition
compatible with the system and involving the condition $\sin
\theta = 0,$ that is, initial condition of the type $(z_{10},
z_{20}, z_{30}, u_{10}, u_{20}, u_{30}, v_{00}) = (0, 0, \pm1,
u_{10}, u_{20}, 0, 0),$ consists of a uniform circular motion of
$z$ on a vertical plane perpendicular to the constant vector
$(u_1(t), u_2(t), u_3(t)) = (u_{10}, u_{20}, 0),$ while $v_0(t) =
0.$
This is also consistent with physical reasoning.\\

\noindent We are not going to perform a detailed study of the
dynamics of the symmetric elastic sphere in this paper, although
it is clear that having the kind of explicit solutions that we
have found clearly helps to study the standard dynamical systems
questions about equilibria, linear and nonlinear stability,
bifurcations, and, also control of the system, with potential
applications to robotics. A future work in this direction is being
planed.
\paragraph{Appendix A.}
We will prove that $M^1_b$ is a nonsingular $3$-dimensional
manifold. We must prove that the matrix
\begin{equation}
\label{matrix} \left[\begin{array}
{ccccccc} 0 & 0 & -v_0 & 0 & 0 & 1 & -z_3 \\
                               0 & 0 & 0 & u_1 & u_2 & u_3 & \lambda v_0 \\
                               z_1 & z_2 & z_3 & 0 & 0 & 0 & 0 \\
                               u_1 & u_2 & u_3 & z_1 & z_2 & z_3 & 0
                               \end{array}\right]\end{equation}
has rank $4$, for all $(z,u,v_0) \in M^1_b.$\\

\begin{enumerate}
\item If {\bf $z_3=\pm 1$}, $z_1=z_2=u_3=v_0=0$, $(u_1,u_2)
\not=0,$ we get
\[
\left[\begin{array}{ccccccc}
0 & 0 & 0 & 0 & 0 & 1 & \pm 1 \\
0 & 0 & 0 & u_1 & u_2 & 0 & 0 \\
0 & 0 & \mp 1 & 0 & 0 & 0 & 0 \\
u_1 & u_2 & 0 & 0  & 0 & \mp 1 & 0 \end{array}\right],\]
which has rank $4.$\\

\item If {\bf $z_3=0$}, $u_3=0$, $(z_1,z_2) \not=0$,
$z_1u_1+z_2u_2=0,$ we have
\[
\left[\begin{array}{ccccccc}
0 & 0 & -v_0 & 0 & 0 & 1 & 0 \\
0 & 0 & 0 & u_1 & u_2 & 0 & \lambda v_0 \\
z_1 & z_2 & 0 & 0 & 0 & 0 & 0 \\
u_1 & u_2 & 0 & z_1  & z_2 & 0 & 0 \end{array}\right].\]

\begin{enumerate}

\item
 If $v_0=0$, then $(u_1,u_2) \not=0.$ For instance, if $u_1 \not=0$
 we have the subdeterminant
\[\left|\begin{array}{cccc}
0 & 0 & 0 & 1 \\
0 & 0 & u_1 & 0 \\
z_1 & z_2 & 0 & 0 \\
u_1 & u_2 & z_1 & 0
\end{array}\right|
=
-u_1\left|\begin{array}{cc}z_1 & z_2 \\
u_1 & u_2 \end{array} \right| \not=0.\] In fact, if
\[
\left|\begin{array}{cc}z_1 & z_2 \\
                            u_1 & u_2 \end{array} \right|=0,
\]
then
\begin{eqnarray*}
-z_2u_1+z_1u_2 & = & 0 \\
z_1u_1+z_2u_2 & = & 0
\end{eqnarray*}
and since $z_1^2+z_2^2=1$, we have $(u_1,u_2)=0$, a contradiction.
A similar result holds if $u_2 \not=0.$
\\

\item
 If $v_0 \not=0$, and, say $z_1 \not=0$, we have
\[\left|\begin{array}{cccc} 0 & -v_0 & 0 & 0 \\
                            0 & 0 & u_1 & \lambda v_0\\
                            z_1 & 0 & 0 & 0 \\
                            u_1 & 0 & z_1 & 0 \end{array} \right|=\lambda v_0\left|\begin{array}{ccc} 0 & -v_0 & 0 \\
                            z_1 & 0 & 0 \\
                            u_1 & 0 & z_1 \end{array} \right|=-\lambda v_0^2z_1^2 \not=0.\]

A similar result holds if $z_2 \not = 0.$\\
\end{enumerate}

\item Assume now that
{\bf $0 < |z_3| < 1$}; $(z_1,z_2)\not=0.$\\

\begin{enumerate}

\item
 If  {\bf $v_0=0;$} $u_3=0,$ $(u_1,u_2) \not = 0$, $z_1u_1+z_2u_2 = 0,$ we have
\[\left[\begin{array}{ccccccc} 0 & 0 & 0 & 0 & 0 & 1 & -z_3 \\
                               0 & 0 & 0 & u_1 & u_2 & 0 & 0 \\
                               z_1 & z_2 & z_3 & 0 & 0 & 0 & 0 \\
                               u_1 & u_2 & 0 & z_1 & z_2 & z_3 & 0\end{array}\right].
                                \]
If, for instance, $u_1 \not=0,$ we have
\[ \left|\begin{array}{cccc} 0 & 0 & 0 & -z_3 \\
                             0 & 0 & u_1 & 0 \\
                             z_1 & z_2 & 0 & 0 \\
                             u_1 & u_2 & z_1 & 0 \end{array}\right|= z_3u_1
\left|\begin{array}
                             {cc} z_1 & z_2\\
                             u_1 & u_2 \end{array}\right|\not=0;\]

In fact, if

$$
\left|\begin{array}{cc} z_1 & z_2\\
                             u_1 & u_2 \end{array}\right|=0$$

then $-z_2u_1+z_1u_2=0$ and this together with $z_1u_1+z_2u_2=0 $
gives $(u_1,u_2)=0$ or
$(z_1,z_2)=0$, a contradiction. Similarly if $u_2 \not=0.$\\

\item If $v_0 \not= 0$; $u_3 \not= 0$, $z_1u_1+z_2u_2 \not= 0$,
$(u_1,u_2) \not=0$; in this case the matrix (\ref{matrix}) also
has rank $4.$ In fact, $(z_1,z_2,z_3)$ and $(u_1,u_2,u_3)$ must be
linearly independent, otherwise the condition
$z_1u_1+z_2u_2+z_3u_3=0$ would imply $u=0.$
Then the last
two rows of (\ref{matrix}) are linearly independent. Using this we
can see easily that the second row cannot be a linear combination
of the last two rows, so the last three rows are linearly
independent.
Let
$R_i,$
$i = 1,2,3,4$
represent the four rows.
If the rank of the matrix is not 4 then
we must have
$R_1 = \chi R_4 + \psi R_3 + \varphi R_2.$
This leads to the systems of equations
\begin{align*}
\chi u_1 + \psi z_1
&= 0\\
\chi u_2 + \psi z_2
&= 0\\
\chi u_3 + \psi z_3
&= -v_0
\end{align*}
and
\begin{align*}
\chi z_1 + \varphi u_1
&= 0\\
\chi z_2 + \varphi u_2
&= 0\\
\chi z_3 + \varphi u_3
&= 1.
\end{align*}
Using
$<z, u> = 0$
we can deduce from the first system that
$\chi u^2 = -v_0 u_3$
and since
$u_3 = v_0 z_3$
we obtain
$\chi u^2 = -v_0^2 z_3.$
On the other hand, using
$<z, u> = 0$
we can deduce from the second system that
$\chi = z_3.$
We can conclude that
$u^2 = - v_0^2,$
which implies
$u = v_0 = 0,$
a contradiction.

Therefore, also
in this case the rank is $4.$
\end{enumerate}
\end{enumerate}

\noindent We can conclude that $M^1_b$ is a nonsingular manifold of
dimension $3,$ defined regularly in $\Bbb R^{7}$ by the equations
(\ref{EEc5})-(\ref{EEc8}).

\paragraph{Appendix B.}
The proof that $f$ is a diffeomorphism goes as follows.

\begin{enumerate}
\item Injectivity

\begin{enumerate}
\item

$z_3 \not= \pm 1.$\\
For given
$(z, (\cos \psi, \sin \psi)) \in S^2 \times S^1,$
we have
$(z,u,v_0)=f(z,(\cos{\psi},\sin{\psi})).$
We have that
$(\theta,\varphi) \in (0,\pi)\times[0,2\pi)$
is uniquely determined, moreover,
$0 < \sin{\theta}$
and
$-1 < \cos{\theta} <  1$
are determined. Therefore from equation (\ref{par7})
$\cos{(\varphi-\psi)}$ is determined. Using (\ref{par4}) and
(\ref{par5}) we obtain $\sin(\varphi-\psi).$ Since $z$ determines
$(\cos{\varphi},\sin{\varphi})$, we get $\psi \in [0,2\pi)$
uniquely determined.

\item
$z_3 = \pm 1.$ \\
In this case $\sin{\theta}=0$, $\cos{\theta}=\pm{1}$,
$a=b=\sqrt{\mu}.$ Using (\ref{par4}) and (\ref{par5}) we obtain
\begin{eqnarray*}
u_1 &=-\sqrt{\mu}\cos \psi\\
u_2 &= -\sqrt{\mu}\sin \psi
\end{eqnarray*}
which determines $\psi \in [0,\,2\pi).$
\end{enumerate}

\item We must prove now that the rank of the tangent map
$T_{(z, (\cos\psi, \sin\psi))}f$
is
$3$
for all
$(z, (\cos\psi, \sin\psi)) \in S^2 \times S^1.$
We have that the Jacobian matrix calculated at
$(\theta, \varphi, \psi)$
of the map given by
(\ref{par1})-(\ref{par7}) is
\[J(\theta, \varphi, \psi)=\left[\begin{array}{ccc}
\cos{\theta}\cos{\varphi} & -\sin{\theta}\sin{\varphi} & 0 \\
\cos{\theta}\sin{\varphi} & \sin{\theta}\cos{\varphi} & 0 \\
-\sin{\theta} & 0 & 0 \\
J_{41} & J_{42}  & J_{43} \\
J_{51} & J_{52} & J_{53} \\
J_{61} & J_{62} & J_{63}\\
a'\cos(\varphi-\psi)\sin{\theta}+a\cos(\varphi-\psi)\cos{\theta} &
-a\sin(\varphi-\psi)\sin{\theta} & a\sin(\varphi-\psi)\sin{\theta}
\end{array}\right]
\]

where we have written
\[\begin{array}{lll}
J_{41} & = & -a'\cos(\varphi-\psi)\cos^2{\theta}\cos{\varphi}
+2a\cos(\varphi-\psi)\cos{\theta}\sin{\theta}\cos{\varphi}\\
J_{42} & = & a\sin(\varphi-\psi)\cos^2{\theta}\cos{\varphi}
+a\cos(\varphi-\psi)\cos^2{\theta}\sin{\varphi}-\\
& &
-b\cos(\varphi-\psi)\sin{\varphi}-b\sin(\varphi-\psi)\cos{\varphi}\\
J_{43}& = & -a\sin(\varphi-\psi)\cos^2{\theta}\cos{\varphi}
+b\cos(\varphi-\psi)\sin{\varphi}\\
J_{51}& = & -a'\cos(\varphi-\psi)\cos^2{\theta}\sin{\varphi}+
2a\cos(\varphi-\psi)\cos{\theta}\sin{\theta}\sin{\varphi}\\
J_{52}& = & a\sin(\varphi-\psi)\cos^2{\theta}\sin{\varphi}
-a\cos(\varphi-\psi)\cos^2{\theta}\cos{\varphi}+\\
& & + b\cos(\varphi-\psi)\cos{\varphi}-b\sin(\varphi-\psi)\sin{\varphi}\\
J_{53} &   = & -a\sin(\varphi-\psi)\cos^2{\theta}\sin{\varphi}-\\
& & -b\cos(\varphi-\psi)\cos{\varphi}\\
J_{61}& = & a'\cos(\varphi-\psi)\cos{\theta}\sin{\theta}
+a\cos(\varphi-\psi)[-\sin^2{\theta}+\cos^2{\theta}] \\
J_{62}& = & -a\sin(\varphi-\psi)\cos{\theta}\sin{\theta}\\
J_{63}& = & a\sin(\varphi-\psi)\cos{\theta}\sin{\theta}
\end{array},\]
with $a'=da /d\theta.$

\begin{enumerate}
\item If $z_3=\pm 1$
then
$\theta = 0$
or
$\theta = \pi.$
Let us consider the case $\theta = 0$ only since the case
$\theta = \pi$
can be treated in a similar way.
We cannot use the polar coordinates representations, in fact, we can see that the rank of
$J(0, \varphi, \psi)$
and
$J(\pi, \varphi, \psi)$
is less than
$3$
for all
$(\varphi, \psi).$
On the other hand,
we can easily check that for fixed
$\psi,$
$f(\theta(t),\varphi,\psi)$ represent the same
point in the manifold $M_b^1$ for $\theta(t)=0$ and
$\varphi$ arbitrary, which is  the point
$(0,0, 1,-b\cos \psi, -b\sin \psi,0,0).$
Let, for instance,
$\theta(t)=t,$ then the tangent
vector to the curve $f(\theta(t),\varphi, \psi)$
at
$t = 0,$
is
$\left(\cos \varphi, \sin \varphi, 0, 0, 0, b\cos(\varphi - \psi), b\cos(\varphi - \psi)\right).$
For
$\varphi = 0$
and
$\varphi = \pi / 2$
we obtain the linearly independent vectors.
$A =\left(1, 0, 0, 0, 0, b\cos \psi, b\cos \psi\right)$
and
$B =\left(0, 1, 0, 0, 0, b\sin \psi, b\sin \psi\right),$
respectively.
On the other hand, the tangent vector to the curve
$f(0, 0, \psi + t)$
at
$t = 0,$
which is the last column of
$J(0, 0, \psi),$
is
$ C = (0, 0, 0, b\sin \psi, -b\cos \psi, 0, 0).$
The three vectors
$A, B, C$
are linearly independent for all
$\psi,$
which shows that the rank of
$T_{(e_3, (\cos \psi, \sin \psi))}f$
is
$3$
for all
$\psi.$

\item If $z_3=0$, $\cos \theta=0$, then $\theta=\pi/2$ and
$a=\sqrt {\mu /\lambda}.$
We get in this case
\[J(\pi/ 2, \varphi, \psi)
=
\left[\begin{array}{ccc}
0 & -\sin{\varphi} & 0 \\
0 &  \cos{\varphi} & 0 \\
-1 & 0 & 0 \\
0 & -b\sin(2\varphi-\psi) & b\cos(\varphi-\psi)\sin \varphi \\
0 &
b\cos(2\varphi-\psi)
& -b\cos(\varphi-\psi)\cos \varphi \\
-a\cos(\varphi-\psi) & 0 & 0 \\
0 & -a\sin(\varphi-\psi) & a\sin(\varphi-\psi)
\end{array}\right] \]

\begin{enumerate}
\item
If $\sin{\varphi}=0$, then $\cos{\varphi}=\pm{1}.$\\
Case I: $\sin(\varphi-\psi) \not = 0$, then the rank of $J$ is $3.$\\
Case II: $\sin(\varphi-\psi) = 0$, then $\cos(\varphi-\psi)=\pm 1$
and $-b\cos(\varphi-\psi)\cos{\varphi}\not= 0$; then the rank of $J$ is $3.$
\item
If $\sin{\varphi}=\pm{1}$, then $\cos{\varphi}=0.$ \\
Case I: $\sin(\varphi-\psi) \not = 0$, then the rank of $J$ is $3.$\\
Case II: $\sin(\varphi-\psi) = 0$, then $\cos(\varphi-\psi)=\pm 1$
and $-b\cos(\varphi-\psi)\sin{\varphi}\not= 0$; then the rank of $J$ is $3.$
\item
If $0 < |\sin{\varphi}| < 1$,\\
Case I:
$\sin(\varphi-\psi) \not = 0$, then the rank of
$J$ is $3.$\\
Case II: $\sin(\varphi-\psi) = 0$, then $\cos(\varphi-\psi)=\pm 1$
and $b\cos(\varphi-\psi)\sin{\varphi}\not= 0$; then the rank of $J$ is $3.$
\end{enumerate}

\item If $0 < |z_3|< 1$, then we have
\[\left|\begin{array}{cc} J_{11} & J_{12} \\
                          J_{21} & J_{22}\end{array}\right|=
                          \left|\begin{array}{cc}\cos{\theta}\cos{\varphi}
                          & -\sin{\theta}\sin{\varphi} \\
 \cos{\theta}\sin{\varphi}& \sin{\theta}\cos{\varphi} \end{array}\right|
 =\cos{\theta}\sin{\theta}\not=0;\]
then the rank of J will be $3$ if $J_{i3}\not=0$,
for at least one $i$, $4 \leq i \leq 7$.\\

\begin{enumerate}
\item
If $J_{73}\not=0$, then the rank of $J$ is $3.$\\

\item If $J_{73}=0$, then $\sin(\varphi-\psi)=0$ and
$\cos(\varphi-\psi)=\pm 1.$ Moreover, $J_{43}=\pm b\sin{\varphi}$
and
$J_{53}=\pm a\cos^2{\theta}\sin{\varphi}\pm b\cos{\varphi}.$\\
Case I: If $\sin{\varphi} \not= 0$, then the rank of $J$ is $3.$\\
Case II: If $\sin{\varphi}=0$, $\cos{\varphi}=\pm 1$ and
$J_{53}=\pm b \not= 0,$ then the rank of $J$ is $3.$
\end{enumerate}
\end{enumerate}
\end{enumerate}
\paragraph{Acknowledgements} We thank Jerrold Marsden and Tudor Ratiu for helpful
comments and references. The first author also thanks the
Bernoulli Center, at EPFL, for its kind hospitality during the
last period of the realization of this work and also the Universidad Carlos III de Madrid for its kind hospitality during part of the realization of this work.


\begin{thebibliography}{Dillon 83}

\bibitem{FM} Abraham, R. and J.E. Marsden, Foundations of
Mechanics, Benjamin, 1978.

\bibitem{MS}Marsden, J.E. and T.S. Ratiu, Introduction to Mechanics
and Symmetry, Springer-Verlag, 1994.

\bibitem{CEMARA1} H. Cendra, J.E. Marsden, T.S. Ratiu, Lagrangian Reduction
by Stages, Memoirs of the AMS 152 (722) (2001).

\bibitem{BKMM} A.M. Bloch, P.S. Krishnaprasad, J.E. Marsden, R.M. Murray,
Nonholonomic Mechanical Systems with Symmetry, Arc. Rat. Mech.
Anal. 136 (1996) 21-99.

\bibitem{CEMARA2} H. Cendra, J.E. Marsden, T.S. Ratiu, Geometric Mechanics,
Lagrangian Reduction and Nonholonomic Systems, in: Mathematics
Unlimited-2001 and Beyond, Springer, pp. 221-273.

\bibitem{CORTES} J. Cortes Monforte, Geometric, Control and
Numerical Aspects of Nonholonomic Systems, Sringer, 2002.

\bibitem{DIRAC} P.A.M. Dirac, Generalized Hamiltonian Dynamics,
Can. J. of Math. 2 (1950) 129-148.

\bibitem{new:the} R. Newcomb, The semistate description of nonlinear
and time-variable circuits, IEEE Trans. Circ. Syst. 28 (1981)
62-71.

\bibitem{zab:pow} V. Venkatasubramanian, H. Schattler, J. Zaborsky,
Local bifurcations and feasibility regions in
differential-algebraic systems, IEEE Trans. Autom. Control 40 (12)
(1995) 1992-2013.

\bibitem{McC:sps} N. McClamroch, H. Krishnan, Nonstandard singularly
perturbed control systems and differential-algebraic equations,
Int. J. Cont. 55 (4) (1992) 1239-1253.

\bibitem{kum:dis} A. Kumar and P. Daoutidis, A DAE Framework for
Modelling and Control of Reactive Columns, 4th IFAC Symp. on Dyn.
Cont. Chem. Reactors, Helsingor, Denmark (1995) 99-104.

\bibitem{kri:robot} H. Krischnan and H. McClamroch, Tracking in
Nonlinear Differential-Algebraic Control Systems with Applications
to Constrained Robot Systems, Automatica 30 (1994) 1885-1897.

\bibitem{rhe2:Impdif} P. Rabier and W. Rheinboldt,
Treatment of Implicit differential-algebraic Equations, J.
Differential Equations 109 (1994) 110-146.

\bibitem{rhe1:difalg} W. Rheinboldt, Differential algebraic systems as
differential equations on manifolds, Math. of Comp. 43 (1984)
473-482.

\bibitem{re1:onag} S. Reich, On a geometrical interpretation of
differential-algebraic equations, Circ. Syst. Sign. Proc. 9 (1990)
367-382.

\bibitem{re2:onae} S. Reich, On an existence and uniqueness theory for
nonlinear differential algebraic equations, Circ. Syst. Sign.
Proc. 10 (1991) 343-359.

\bibitem{rab:gen} P. Rabier and W. Rheinboldt,
A general existence and uniqueness theory for implicit
differential algebraic equations, Diff. Int. Equations. 4 (1991)
563-582.

\bibitem{sza:geo} A. Szatkowski, Geometric Characterization of
Singular Differential Algebraic Equations, Int. J. Syst. Sci. 23
(1992) 167-186.

\bibitem{pritch:onimp} F. Leon Pritchard, On implicit
systems of differential equations, J. Differential Equations 194
(2003) 328-263.

\bibitem{rei:sin} G. Reibig, On Singularities of
Autonomous Implicit Ordinary Differential Equations, IEEE Trans.
Circ. Syst. 50 (2003) 922-931.

\bibitem{soto:imp} J. Sotomayor, Impasse Singularities of Differential
Systems of the Form $A(x)\,x'=F(x)$, J. Differential Equations 169
(2001) 567-587.

\bibitem{guz:imp} A. Guzman Gomez, Constrained Equations with
Impasse Points, J. Math. Anal. Appl. 214 (1997) 292-306.

\bibitem{rhe3:imp} P. Rabier, W. Rheinboldt,
On the computation of impasse points of quasilinear
differential-algebraic equations, Math. Comp. 62 (1994) 133-154.

\bibitem{hiron1:desing} H. Hironaka, Resolution of an
algebraic variety over a field of characteristic zero, Ann. of
Math. 79 (1964) 109-203.

\bibitem{hau:Hironaka} H. Hauser,
The Hironaka Theorem on Resolution of Singularities, Bulletin of
the AMS 40 (2003) 323-403.

\bibitem{villa:new} S. Encinas and O. Villamayor,
A new proof of desingularization over fields of characteristic
zero, Rev. Mat. Iber. 19 (2003) 339-353.

\bibitem{BIMI} E. Bierstone, P. D. Milman,
Canonical desingularization in characteristic zero by blowing up
the maximum strata of a local invariant, Invent. Math 128 (1997)
207-302.

\bibitem{BM2003} E. Bierstone, P. D. Milman,
Desingularization algorithms I. Role of exceptional divisors.
Moscow Math. Journal, Volume 3 (3) (2003) 751-805.

\bibitem{hiron2:subanal} H. Hironaka, Subanalytic Sets,
Number Theory, Algebraic Geometry and Commutative Algebra, Tokio,
Kinokuniya (1973) 453-493.

\bibitem{bier88:anal} E. Bierstone and P.D. Milman,
Semianalytic and subanalytic sets, Publications mathem\'{a}tiques
de l'I.H.E.S., tome 67 (1988) 5-42.

\bibitem{suss98:control} H.J. Sussmann,
Some optimal control applications of real-analytic stratifications
and desingularization, Singularities Symposium-Lojasiewicz70,
Banach Center Publ. 44 (1998), Polish. Acad. Sci., Warsaw.


\bibitem{suss90:desing} H.J. Sussmann, Real analytic desingularization
and subanalytic sets: an elementary approach, Trans. Amer. Math.
Soc. 317 (2) (1990) 417-461.

\bibitem{vesel:newcases} L.E. Veselova, New cases of the integrability of the
equations of a rigid body in the presence of a nonholonomic
constraint, Geometry, differential equations and mechanics, Moscow
(1985) 64-68.

\bibitem{MTA} Abraham, R., J.E. Marsden and T.S. Ratiu,
Manifolds, Tensor Analysis and Applications, Springer-Verlag,
1988.

\bibitem{rabier:lineal} P. Rabier, W. Rheinboldt,
Classical and generalized solutions of time-dependent linear
DAEs., Lin. Alg. Appl. 245 (1996) 259-293.

\bibitem{lojasiewicz1} S. Lojasiewicz, Sur le Probl\'{e}me de la division,
Studia Math.18 (1959) 87-136.

\bibitem{lojasiewicz2} S. Lojasiewicz, Triangulation of Semianalytic sets,
Ann. Scuola. Norm. Sup. Pisa 18, (1964) 449-474.

\bibitem{lojasiewicz3} S. Lojasiewicz, Sur les ensembles semi-analytiques,
Actes du Congress International des Mathematiciens, Nice (1970),
Tome 2, 237-241.

\bibitem{Gabrielov} A. M. Gabrielov, Formal relations between analytic
functions, Math. USSR Izvestija 7 (1973) 1056-1088.

\bibitem{Hardt1} R.M. Hardt, Stratifications of real analytic mappings
and images, Inven. Mathem. 28 (1975) 193-208.

\bibitem{Hardt2} R.M. Hardt, Triangulation of subanalytic sets and proper
light subanalytic maps, Inven. Mathem. 38 (1977) 207-217.

\bibitem{sussmann73} H. Sussmann, Orbits of families of vector fields
and integrability of distributions, Trans. Amer. Math. Soc. 180
(1973) 171-188.

\bibitem{fufaev} Neimark, Ju. I. and N. A. Fufaev,
Dynamics of Nonholonomic Systems, Translations of Mathematical
Monographs 33 AMS, Providence RI (1972).

\bibitem{CMDLI} H. Cendra, A. Ibort, M. de Le{\'o}n, D. Martin de Diego,
A generalized Chetaev's principle for a class of higher order
nonholonomic constraints, J. Math. Phys. 45 (7) (2004) 2785-2801.

\bibitem{CELARE} H. Cendra, E.A. Lacomba, W.A. Reartes,
The Lagrange D'Alembert-Poincar\'{e} equations for the symmetric
rolling sphere, Actas del VI Congreso Antonio, A.R Monteiro,
Universidad Nacional del Sur, Bah{\'\i}a Blanca (2002) 19-32.
\end{thebibliography}
\end{document}